\documentclass[11pt]{amsart}
\usepackage{adjustbox, amssymb, amsfonts, amsmath, bbm, bm, mathrsfs, mathtools, pgfplots, physics, stackrel, tikz, tikz-cd,  todonotes,  verbatim, xypic}
\usetikzlibrary{intersections, calc, matrix, decorations.pathreplacing}

\usepackage[shortlabels]{enumitem}
\usepackage[bbgreekl]{mathbbol}
\DeclareMathSymbol\bbDelta  \mathord{bbold}{"01}

\usepackage[outline]{contour} \contourlength{2pt}
\usepackage{hyperref}
\hypersetup{
    colorlinks=false, %set true if you want colored links
    linktoc=all,     %set to all if you want both sections and subsections linked
    linkcolor=blue,  %choose some color if you want links to stand out
}

\pgfplotsset{compat=newest}

\usepackage{geometry}
\allowdisplaybreaks

\theoremstyle{definition}
\newtheorem{theorem}{Theorem}[section]
\newtheorem{lemma}[theorem]{Lemma}
\newtheorem{corollary}[theorem]{Corollary}
\newtheorem{proposition}[theorem]{Proposition}
\newtheorem{definition}[theorem]{Definition}

\newtheorem{example}[theorem]{Example}
\newtheorem{notation}[theorem]{Notation}
\newtheorem{remark}[theorem]{Remark}
\numberwithin{equation}{section}
\numberwithin{figure}{section}

\numberwithin{equation}{section}

\DeclareMathOperator*{\sSet}{{\it sSet}}

\DeclareMathOperator*{\Tot}{Tot}
\DeclareMathOperator*{\tot}{tot}

\DeclareMathOperator*{\sPre}{Pre^{\infty}}
\DeclareMathOperator*{\modmod}{/ \! \!/}
\newcommand{\mmod}[2]{{#1}\!  \modmod  \! {#2}}

\newcommand\subsetsim{\mathrel{%
  \ooalign{\raise0.2ex\hbox{$\subset$}\cr\hidewidth\raise-0.8ex\hbox{\scalebox{0.9}{$\sim$}}\hidewidth\cr}}}
\newcommand{\nGpd}{\Gamma(\mathbb{C}^n)}
\newcommand{\Bihol}{\text{St}\mathbb{C}^n}

\newcommand{\Chart}{ \mathcal{C}hart}

\newcommand{\del}{ \partial}
\newcommand{\delbar}{ \overline{\partial}}

\newcommand{\Cech}{\check{C}}

\newcommand{\extAnti}{\mathbb{\Omega}}
\newcommand{\holForm}{\Omega_{\text{hol}}}
\newcommand{\BMmap}{\mathbf{BM}}
\newcommand{\Hartmap}{\mathbf{Hart}}

\newcommand{\holFormP}[1]{{\it\Omega^{#1}_{\text{\it hol}}}}
\newcommand{\DK}{{DK}}
\newcommand{\Z}{\mathbb Z}
\newcommand{\C}{\mathbb C}
\newcommand{\T}{\mathcal T}
\newcommand{\bu}{{\bullet}}
\newcommand{\ot}{\otimes}
\newcommand{\wt}{\widetilde}
\newcommand{\mc}{\mathcal}
\newcommand{\Cfmap}[1]{\mathbf{C}_{#1}}

\newcommand{\Om}{\Omega}

\newcommand{\OAcpt}{A_{\text{cpt}}}
\newcommand{\om}{\omega}
\newcommand{\al}{\alpha}
\newcommand{\be}{\beta}
\newcommand{\Del}{\Delta}

\newcommand{\la}{\langle}
\newcommand{\ra}{\rangle}

\DeclareFontFamily{U}{dmjhira}{}
\DeclareFontShape{U}{dmjhira}{m}{n}{
  <-> s*[0.95] dmjhira
}{}
\DeclareFontSubstitution{U}{dmjhira}{m}{n}

\DeclareFontFamily{OT1}{pzc}{}
\DeclareFontShape{OT1}{pzc}{m}{it}{<-> s * [1.10] pzcmi7t}{}
\DeclareMathAlphabet{\mathpzc}{OT1}{pzc}{m}{it}

\newcommand{\yo}{\mathpzc{y}}
\newcommand\restr[2]{{% we make the whole thing an ordinary symbol
  \left.\kern-\nulldelimiterspace % automatically resize the bar with \right
  #1 % the function
  \vphantom{\big|} % pretend it's a little taller at normal size
  \right|_{#2} % this is the delimiter
  }}

\author[C. Glass]{Cheyne Glass}
\address{Cheyne Glass, State University of New York at New Paltz, Department of Mathematics, 1 Hawk Dr., New Paltz, NY 12561}
  \email{glassc@newpaltz.edu}

\author[T. Tradler]{Thomas Tradler}
  \address{Thomas Tradler, New York City College of Technology The City University of New York, Department of Mathematics, 300 Jay Street, Brooklyn, NY 11201}
  \email{ttradler@citytech.cuny.edu}

\author[M. Zeinalian]{Mahmoud Zeinalian}
  \address{Mahmoud Zeinalian, Lehman College, The City University of New York, Department of Mathematics, 250 Bedford Park Blvd W, Bronx, NY 10468}
  \email{mahmoud.zeinalian@lehman.cuny.edu}

\title[Hirzebruch-Riemann-Roch for complex analytic $\infty$-prestacks]{Hirzebruch-Riemann-Roch for\\ complex analytic $\infty$-prestacks}

\begin{document}
\maketitle
\begin{center} 
\vspace{-4mm}{\it Dedicated to Domingo Toledo and Yue Lin Lawrence Tong}
\end{center}

\vspace{1mm}
\begin{abstract}
We provide a cocycle-level Hirzebruch-Riemann-Roch (HRR) identity for arbitrary complex analytic $\infty$-prestacks. We view this work as the natural setting for Toledo and Tong's HRR philosophy and technical machinery.  
\end{abstract}

\setcounter{tocdepth}{1}
\tableofcontents
\section*{Introduction}
The celebrated Hirzebruch-Riemann-Roch theorem (HRR) gives a formula for the Euler characteristic or the alternating sum of the ranks of the $\delbar$ cohomology of a holomorphic vector bundle on a complex manifold, in terms of the integral of the Todd class of the base times the Chern character of the bundle. The alternating sum of the dimensions of the cohomology groups agrees with the alternating sum of the dimensions of the terms in the $\delbar$ resolution of holomorphic sections of the bundle. 

While this resolution is infinite-dimensional in each degree, there is a precise cancellation of infinities, resulting in the desired finite answer. One way to make this precise is to think of the dimension as the trace of the identity map and deform the identity map up to homotopy to a smoothing operator that is trace class. Such a homotopy between the identity map and a smoothing operator on the $\delbar$ complex of forms with values in a holomorphic vector bundle is called a $\delbar$ parametrix, or simply a parametrix.

The kernel of the smoothing operator is a form that lives on the Cartesian product of the manifold with itself, and standard Chern-Weil theory of invariant polynomials identifies the restriction of the kernel to the diagonal with the Todd class times the Chern character of the bundle. When the manifold is compact, integration gives the desired HRR as the equality of two integers.

The most common tool for finding a parametrix for a differential operator is to choose an auxiliary metric and use its heat flow to deform the identity map to a smoothing operator. This elegant technique has inherent limitations. For instance, when a discrete group acts on a manifold, one wants to treat everything equivariantly. Unless the group is finite, there is no hope of finding an invariant metric on the manifold. We use alternative techniques developed by Toledo and Tong \cite{TT76} in the 1970s to obviate such difficulties and extend HRR far beyond the equivariant setting.

In this paper, we introduce a site $\Bihol$, whose objects are Stein manifolds biholomorphic with an open subset of $\mathbb{C}^n$, and whose morphisms are holomorphic embeddings. There is a natural site structure (as in definition \ref{DEF site}) given by families $\{U_i \to U\}_{i \in I}$ of morphisms in $\Bihol$ whose images cover $U$. A natural $\infty$-stack, $\Chart$, is defined on this site, formalizing the idea of local coordinates and their transition functions. We repurpose Toledo-Tong's construction of a \v{C}ech parametrix, extending the Bochner-Martinelli kernel, as a map of $\infty$-prestacks \[\BMmap: \Chart \to \extAnti,\] where $\extAnti$ is an $\infty$-prestack built from certain differential forms defined near the diagonal in the Cartesian square of an object of $\Bihol$. Assuming $n\ge 2$, the Hartogs extension theorem yields another map of $\infty$-prestacks  \[\Hartmap: \extAnti \to \holFormP{n},\] where $\holFormP{n}$ is an $\infty$-stack sending an object in $\Bihol$ to its holomorphic $n$-forms, regarded as a simplicial set via the Dold-Kan correspondence. Guided once again by Toledo-Tong, we adapt the Chern-Weil mechanism to define a third map \[ \Cfmap{Todd_n} : \Chart \to  \holFormP{n}.\]

Our Hirzebruch-Riemann-Roch theorem is then the statement that the two maps $\Hartmap \circ \BMmap$ and $\Cfmap{Todd_n}$ are equal as maps of $\infty$-stacks. In this way, we obtain an HRR theorem that applies very generally: for any $\infty$-prestack $\mc F$ on $\Bihol$, HRR is expressed as an equality of two maps of spaces (homotopy types)
\begin{equation*}
\begin{tikzcd}[column sep=huge]
\mathbb{R}\text{Hom}( \mc F, \Chart)
  \arrow[bend left=15]{r}[name=U, label=above:$\left(\Hartmap \circ \BMmap\right)_*$]{}
  \arrow[bend right=15]{r}[name=D,label=below:${\Cfmap{Todd_n}}_*$]{} &
\mathbb{R}\text{Hom}( \mc F, \holFormP{n})
  \arrow[shorten <=20pt,shorten >=15pt,equal,to path={(U) -- (D)}]{}.
\end{tikzcd}
\end{equation*} We also explain how the induced map on $\pi_0$ recovers familiar classical statements. For example, when $\mc F=\varepsilon M$ is the $\infty$-prestack assigning to each object $W \in \Bihol$ the discrete simplicial set $\varepsilon M(W)$ of holomorphic embeddings $W \to M$, our general HRR theorem asserts the equality of two \v{C}ech $n$-cocycles with values in holomorphic $n$-forms on $M$. Choosing a partition of unity would translate this equality into an equality of two $2n$-forms on $M$. When $M$ is compact, integrating these top forms yields two numbers that are the two sides of the classical HRR theorem,
\begin{equation*}
\int_M td(M)=\chi_{\delbar}(M).
\end{equation*} 
This relationship to the classical HRR theorem in the context of the work by Toledo and Tong is reviewed in more detail below without appealing to the language of $\infty$-prestacks.

For a quotient stack $\mmod{M}{G}$, our setting recovers and generalizes equivariant HRR, when $G$ is a discrete group. While the integration step in this case requires considerations beyond the current work, the HRR takes the form of the cocycle level equality of two classes in the $G$-equivariant cohomology of $M$. When $M= \mathbb{C}^n$, the results of this paper are compatible with the prequel to this paper \cite{GZ}. 

\subsection*{Outline of the paper}
The structure of the paper is as follows. In Chapter 1, the site structure on $\Bihol$ is defined, and the basics of $\infty$-prestacks on $\Bihol$ are introduced. Chapter 2 defines the $\infty$-stack $\Chart$ of local charts. Chapter 3 defines two other $\infty$-prestacks, $\extAnti$ and $\holFormP{n}$, each constructed from a chain complex of differential forms, and then the Hartogs extension theorem is interpreted as a map of $\infty$-prestacks $\Hartmap: \extAnti \to \holFormP{n}$. Here, we use the assumption throughout the paper that $n\ge  2$. In Chapter 4, we interpret a construction of Toledo and Tong based on the Bochner-Martinelli $\delbar$-parametrix as a map of $\infty$-prestacks $\BMmap: \Chart \to \extAnti$. Chapter 5 contains a strong form of the Hirzebruch-Riemann-Roch theorem (theorem \ref{THM main theorem}) for arbitrary $\infty$-stacks over the site $\Bihol$. Chapter 6 provides a specific cofibrant replacement for the quotient $\infty$-prestack of a discrete group acting on a complex manifold. Chapter 7 gives a geometric application in terms of obstructions to conjugating the action of a discrete group on a manifold to an affine one. 

\subsection*{A brief review of Toledo and Tong's HRR machinery}
It is worthwhile outlining Toledo and Tong's proof of the Hirzebruch-Riemann-Roch theorem in \cite{TT76} via a \v{C}ech parametrix, at least for the special case of a trivial $\C$-bundle over $M$ considered in this paper. The sketch given here also serves for much of the motivation of the constructions presented in chapters \ref{SEC presheaves on Bihol}-\ref{SEC:Invariant-polynomials} below, which give a reinterpretation of this in terms of maps of $\infty$-prestacks.

For a given atlas $\{\rho_i:U_i\to \C^n\}_{i\in I}$ of the complex manifold $M$, Toledo and Tong construct a sequence of forms $\om^q=\{\om^q_{i_0,\dots,i_q}\in  A^{{(n,n-q-1)}}(W_{i_0,\dots,i_q}-\Del)\}_{i_0,\dots,i_q\in I}$ for $q=0, \dots, n-1$, where the open sets $W_{i_0,\dots,i_q}\subset \bigcap_{j=0}^q (U_{i_j}\times U_{i_j})$ are certain neighborhoods of the diagonal $\Delta\hookrightarrow M\times M$ contained in $\bigcap_{j=0}^q (U_{i_j}\times U_{i_j})$. Here, the initial forms $\omega^0_{i_0}$ are pullbacks of the Bochner-Martinelli kernel \eqref{DEF Bochner-Martinelli kernel} using the charts $\rho_{i_0}:U_{i_0}\to \C^n$. These forms satisfy $\delbar\om^0_{i_0}=0$, while the higher forms are constructed so that
\begin{equation}\label{EQU:deltawq-1=delbarwq}
\delta\om^{q-1}_{i_0,\dots,i_q}:=\sum_j (-1)^j \om^{q-1}_{i_0,\dots,\widehat{i_j},\dots,i_q}\big|_{W_{i_0,\dots,i_q}}=\delbar \om^q_{i_0,\dots,i_q}\quad\quad \text{ for }q=1,\dots,n-1.
\end{equation}
The sequence $(\om^0,\dots, \om^{n-1})$ is called the \v{C}ech parametrix. Moreover, Toledo and Tong denote by $\delta \om^{n-1}$ by $\om^n$, made up of forms $\om^n_{i_0,\dots, i_n}\in A^{(n,0)}(W_{i_0,\dots, i_n}-\Delta)$, which are holomorphic due to \eqref{EQU:deltawq-1=delbarwq}. In the case $n\geq 2$, the $\om^n_{i_0,\dots, i_n}$ extend holomorphically over the diagonal, which can then be pulled back via the diagonal inclusion $\Delta\hookrightarrow M\times M$. One of the main theorems in \cite{TT76} is the Universal Riemann-Roch theorem \cite[corollary 11.5]{TT76} (and its bundle-valued generalization \cite[corollary 12.1]{TT76}), which asserts that the pullback of the symmetrized $\omega^n$ forms equals the $n^{\text{th}}$ Todd polynomial applied to appropriate Maurer-Cartan forms, i.e., $\Delta^*(\tilde\om^n) =Todd_n(\theta^n)$.

Using a suitable partition of unity, Toledo and Tong assemble the above forms on the given open subsets to a global form via a map $\mathscr D:\prod\limits_{i_0,\dots,i_q}A^{(n,q')}(W_{i_0,\dots,i_q}-\Del)\to A^{(n,q'+q)}(M\times M-\Del)$, for which $\mathscr D\circ (\delbar-\delta)(c)$ equals $(-1)^q\cdot \delbar\circ \mathscr D(c)$ up to a smooth form that vanishes near the diagonal, where $c$ is an element concentrated in $A^{(n,q')}(W_{i_0,\dots,i_q}-\Del)$. Applying $\mathscr D$ to the \v{C}ech parametrix, they define $k=\sum_{q=0}^{n-1}(-1)^q\mathscr D(\tilde\om^q)\in A^{(n,n-1)}(M\times M-\Del)$, and, using this, they also define a smooth form $s\in A^{(n,n)}(M\times M)$ which is given off the diagonal by $s|_{M\times M-\Del}=-\delbar(k)\in A^{(n,n)}(M\times M-\Del)$ and which extends smoothly over the diagonal, since $-\delbar(k)$ equals $\mathscr D(\tilde\om^{n})$ near the diagonal. A crucial identity that holds for these forms $k$ and $s$ is that for all compactly supported forms $\gamma\in\OAcpt^{(n,n-q),(0,q)}(M\times M)$ one has the identity
\begin{equation}\label{EQU:intk=intD-ints}
\int_{M\times M} k\cdot \delbar\gamma=\int_M \Del^*\gamma-\int_{M\times M} s\cdot \gamma,
\end{equation}
where we now also use $\Del:M\to M\times M$ to denote the diagonal map. Considering the associated operators $K$ and $S$ with kernels $k$ and $s$ respectively, i.e., $S:\OAcpt^{(0,q)}(M)\to A^{(0,q)}(M)$, $S(\alpha)(x)=\int_{y\in M} s(x,y)\cdot \alpha(y)$ and $K:\OAcpt^{(0,q)}(M)\to A^{(0,q-1)}(M)$, $K(\alpha)(x)=\int_{y\in M} k(x,y)\cdot \alpha(y)$, equation \eqref{EQU:intk=intD-ints} can be shown to be equivalent to $\delbar K+K \delbar=1-S$, where $S$ is a smoothing operator. Because of this, $K$ is also called a parametrix, now given as an operator on forms. Moreover, from the Universal Riemann-Roch theorem it follows that the class of $\Delta^* s$ coincides with the degree $2n$-part of the Todd class $td(M)$ of $M$.

Now, for a compact complex manifold $M$, recall the non-degenerate Serre duality pairing $\la-,- \ra : H_{\delbar}^{n-p,n-q}(M)\ot H_{\delbar}^{p,q}(M)\to \C$, $\la \al,\be\ra=\int_M \al\cdot \be$ for Dolbeault cohomology. We note that, in this terminology, \eqref{EQU:intk=intD-ints} can be interpreted by saying that the class $[s]$ on $M\times M$ is the Serre dual of the diagonal $\Del$ interpreted as a distribution $[\gamma]\mapsto \int_\Del \gamma$. Writing $\chi_{\delbar}(M)=\sum_{q=0}^n (-1)^q \dim H_{\delbar}^{0,q}(M)$ for the stated Euler characteristic, one obtains the usual Hirzebruch-Riemann-Roch theorem for the trivial $\C$-vector bundle over $M$:
\begin{equation}\label{EQU:int_todd=Euler}
\int_M td(M)=\chi_{\delbar}(M)
\end{equation}
\begin{proof}
For a basis $\{[\psi_i]\}_i$ of $\bigoplus_q H_{\delbar}^{0,q}(M)$ with $\psi_i\in \Om^{0,|\psi_i|}(M)$, let $\{[\psi_i^\sharp]\}_i$ be a dual basis of $\bigoplus_q H_{\delbar}^{n,n-q}(M)$ with $\la \psi_i^\sharp,\psi_j\ra=\int_M \psi_i^\sharp\cdot \psi_j =\delta_{i,j}$. Then, the Serre dual of $[\sum_i \psi_i(x)\cdot \psi_i^\sharp(y)]$ in $M\times M$ is the diagonal $\Del$, since $\int_M\Del^*(\psi^\sharp_\ell(x)\cdot\psi_m(y))=\int_M \psi_\ell^\sharp(x)\cdot \psi_m(x) =\delta_{\ell,m}$ and
\[
\int_{M\times M}\Big(\sum_i \psi_i(x)\cdot \psi_i^\sharp(y)\Big)\cdot (\psi^\sharp_\ell(x)\cdot\psi_m(y))=
\sum_i \int_M \psi^\sharp_\ell(x)\cdot \psi_i(x)\cdot \int_M \psi_i^\sharp(y)\cdot \psi_m(y)=\delta_{\ell,m},
\]
and so, $[s(x,y)]$ equals $[\sum_i \psi_i(x)\cdot \psi_i^\sharp(y)]$ in $\bigoplus_q H_{\delbar}^{(0,q),(n,n-q)}(M\times M)$. We then obtain \eqref{EQU:int_todd=Euler}, since $\int_M\Del^*s=\int_M \sum_i \psi_i(x)\cdot \psi_i^\sharp(x)=\sum_i (-1)^{|\psi_i|}=\sum_q (-1)^q \dim H_{\delbar}^{0,q}(M)=\chi_{\delbar}(M)$.
\end{proof}

\subsection*{The case of coefficients in an arbitrary bundle}
For simplicity, this paper treats only the case in which the vector bundle in the classical Hirzebruch-Riemann-Roch theorem is the trivial line bundle. The same construction extends to arbitrary holomorphic vector bundles by adjoining local holomorphic frame data to the coordinate charts and using the bundle-valued universal parametrix constructed by Toledo and Tong. Given the already substantial and highly technical nature of the paper, we have chosen to work in this simplified setting in order to highlight the central features of our construction. We briefly indicate the necessary modifications for a holomorphic vector bundle of rank $k$, following the methods of section 6 of Toledo and Tong.

The first modification is to enlarge the prestack $\Chart$ to the prestack $\Chart \times B(GL(k, \C)^{(-)})$ so that it also records local holomorphic frames of the bundle. The second factor is the classifying prestack associated with holomorphic maps into the general linear group. When evaluated on a complex manifold, the resulting simplicial set has points that record an atlas together with a rank $k$ holomorphic vector bundle described by transition cocycles on the domains of the charts.

The Bochner-Martinelli construction is modified in parallel. In degree zero, the scalar Bochner-Martinelli kernel is replaced by the same kernel multiplied by the identity matrix, so that it becomes a matrix-valued kernel. The next term of the universal parametrix kills the difference between $\omega^0\cdot Id$ and its transform under a pair consisting of a coordinate change and a change of frame. The coordinate change acts by pullback, while the change of frame acts on the matrix values. Away from the diagonal, the frame changes act separately in the two variables, but on the diagonal, this action reduces to conjugation. Completing this procedure yields a matrix-valued $\BMmap$ map. We then apply Hartogs extension and take the trace of the matrix to obtain a map of prestacks from $\Chart \times B(GL(k, \C)^{(-)})$ to $\holFormP{n}$.

On the characteristic class side, the construction of $\Cfmap{Todd_n}$ is modified by replacing the Todd polynomial with the product of the Todd polynomial and the Chern character polynomial. The bundle-valued version of the Toledo-Tong identity then identifies the resulting cocycle with the degree $n$ component of the product of the Todd class of the tangent bundle and the Chern character of the vector bundle. This recovers the usual Hirzebruch-Riemann-Roch formula for an arbitrary holomorphic vector bundle.

\subsection*{Acknowledgements}
We would like to thank Dennis Sullivan, Peter Teichner, Domingo Toledo, Tim Hosgood, Emilio Minichiello, Owen Gwilliam, Brian Williams, Julien Grivaux, Damien Calaque, and Manuel Rivera for helpful conversations over the years regarding the work of Toledo and Tong. We also thank the referee for insightful and detailed recommendations which have improved the exposition of the paper. The second author was partially supported by a PSC-CUNY research grant.

\section{$\infty$-prestacks on $\Bihol$}\label{SEC presheaves on Bihol}
In this section, we study $\infty$-prestacks and $\infty$-stacks on a suitable category of Stein manifolds $\Bihol$. More precisely, $\Bihol$ is defined to be the category whose objects are $n$-dimensional Stein manifolds which are biholomorphic with an open subset of $\mathbb{C}^n$, and whose morphisms are holomorphic embeddings, i.e., biholomorphisms onto their image.

Recall (for example \cite[theorem 11.6.1]{Hi}) that given an (essentially) small category $\mathscr{C}$, the category $\sPre(\mathscr{C})$ of \emph{$\infty$-prestacks} (often referred to as simplicial\footnote{The choice of language here is to eventually avoid confusion between a simplicial object in the category of sheaves, and a simplicial presheaf which is fibrant in the relevant local model structure, both of which might be called ``simplicial sheaves''.} presheaves), $\mc F : \mathscr{C}^{op} \to \sSet$, on $\mathscr{C}$ has a simplicial \emph{(global) projective model structure} whose fibrations are object-wise Kan fibrations of simplicial sets, weak-equivalences are object-wise weak equivalences of simplicial sets (inducing isomorphisms on homotopy groups), and cofibrations are those satisfying a lifting property with respect to the fibrations and weak equivalences. The simplicial enrichment is given by the usual formula, 
\begin{equation}\label{EQ sHom}
\underline{\sPre}(\mathscr{C})\left( \mathcal F, \mathcal G \right)_k:= {\sPre(\mathscr{C})}\left( \mathcal F \otimes \Delta^k, \mathcal G \right),
\end{equation}
utilizing the usual tensor of $\sPre$ over $\sSet$.

In order to discuss sheaves/stacks, a site structure on $\mathscr{C}$ is required.
\begin{definition}\label{DEF site}
A \emph{site} \cite[definition 0.1]{Ar} consists of a category $\mathscr{C}$ together with a collection $\text{Cov}(\mathscr{C})$ of families $\{ U_i \xrightarrow{f_i} U\}_{i \in I}$ of specified morphisms in $\mathscr{C}$ called \emph{covering families}, satisfying:
\begin{enumerate}[(i)]
\item If $f: U' \to U$ is an isomorphism then the singleton set $\{ U' \xrightarrow{f} U\}$ is in $\text{Cov}(\mathscr{C})$.
\item If $\{U_i \to U\}_{i \in I}$ is in $\text{Cov}(\mathscr{C})$ and for each $i \in I$, $\{V_{i,j} \to U_i\}_{j \in J_i}$ is in $\text{Cov}(\mathscr{C})$, then the family $\{V_{i,j} \to U_i \to U\}_{i \in I, j \in J_i}$ is again in $\text{Cov}(\mathscr{C})$. 
\item If $\{U_i \to U\}_{i \in I}$ is in $\text{Cov}(\mathscr{C})$ and $V\to U$ is any morphism in $\mathscr{C}$, then $U_i \times_U V$ exists in $\mathscr{C}$ and the family $\{U_i \times_U V \to V\}_{i \in I}$ is again in $\text{Cov}(\mathscr{C})$.
\end{enumerate}
Call a morphism $V \to U$ in $\mathscr{C}$ \emph{basal} if it belongs to some covering family, then a \emph{Verdier site} \cite[definition 9.1]{DHI} additionally satisfies 
\begin{enumerate}[(i), resume]
\item If $V \to U$ is a basal morphism, then the diagonal $V \to V \times_U V$ is also basal.
\end{enumerate}
\end{definition}
Note that, as usual, $\Bihol$ is assumed to have an initial object so that ``empty pullbacks'' exist. While the category of smooth Cartesian spaces with good covers as covering families does not satisfy axiom (iii) above, the category $\Bihol$ in contrast offers both the restricted morphisms (which do not change which prestacks are stacks) as well as the properties of Stein manifolds, the site structure on $\Bihol$ defined below satisfies all four axioms above. The main tool for this claim is \cite[lemma 4.1]{L} which proves that the inverse image of a Stein open subset of any manifold, under a holomorphic map out of a Stein manifold, is again a Stein open subset. As L\'{a}russon points out, an application of this lemma is that the non-empty intersection of finitely many Stein open subsets is again Stein. Importantly, this condition is exactly the one that fails in the smooth Cartesian setting. Namely, if $\mc U$ and $\mc V$ are two good open covers of a topological space $X$, so that each finite intersection $U_{i_0} \cap \cdots \cap U_{i_p}$ and $V_{j_0} \cap \cdots \cap U_{j_q}$ is always again contractible (i.e., the covers are good), there is no reason (and plenty of counterexamples) why every $U_i \cap V_j$ ought to be contractible. 

\begin{definition}\label{DEF cov fam Bihol}
A family of morphisms $\mathcal{U}=\{U_i \xrightarrow{f_i} U\}_{i \in I}$ in $\Bihol$ is a \emph{covering family} if $\cup_{i \in I} f(U_i)=U$. The collection of all covering families is denoted $\text{Cov}(\Bihol)$.
\end{definition}

\begin{lemma}
The collection $\text{Cov}(\Bihol)$ provides a Verdier site structure on $\Bihol$ .

\begin{proof}
The cited definition for a site has three conditions to check, translated here into the setting of $\Bihol$ for the readers' convenience. The first is that any biholomorphism in $\Bihol$ forms a singleton covering family, which is immediately true. The third is that if $W \xrightarrow{g} U$ is a holomorphic embedding between objects in $\Bihol$ and $\mathcal{U}= \{U_i \xrightarrow{f_i} U\}_{i \in I}$ is a covering family of $U$, then each pullback  $W \times_{U} U_i = g^{-1}\left( f (U_i)\right)$ exists and the induced family $g^*\mathcal{U}:= \{ g^{-1}\left( f (U_i)\right) \xrightarrow{g_i} W \}_{i \in I}$, where each $g_i$ is the restriction of $g$, needs to be a covering family. The fact that $g^*\mc U$ is a cover by objects in $\Bihol$ is in the proof of \cite[lemma 4.1]{L}. The second condition is that if  $\mathcal{U}= \{U_i \xrightarrow{f_i} U\}_{i \in I}$ is a covering family of $U$ and for each $i\in I$  $\mathcal{V}_{i}= \{V_{i,{\alpha}} \xrightarrow{f_{i,\alpha}} U_i\}_{\alpha \in A_i}$ is a covering family for $U_i$, then $\mathcal{V}_{\bullet}=\{V_{i,{\alpha}} \xrightarrow{f_{i,\alpha}} U_i \xrightarrow{f_i} U\}_{i\in I, \alpha \in A_i}$  needs to be a covering family; which follows by the aforementioned comment by L\'{a}russon's after \cite[lemma 4.1]{L}. The fourth condition in order to be a Verdier site, is that if $U_i \to U$ is a map in one of the covering families then $U_i \to U_i \times_U U_i$ is again in one of the covering families. However, due to our morphisms being injections, this map is the identity which is indeed a holomorphic embedding and, by definition of our covering families and the first property, always in some covering family.
\end{proof}
\end{lemma}

Following \cite[theorem 6.2]{DHI} there are two ``local projective model structures'' on $\sPre(\Bihol)$ which actually coincide: the one where the weak equivalences are given by the topological weak equivalences of Jardine (see the proof of our proposition \ref{PROP cofib repl M mmod G} later) and the other presented as a left Bousfield localization at the class $\mathcal{S}$ of hypercovers.

\begin{definition}\label{DEF loc proj mod}
Let $\sPre(\Bihol)$ be the simplicial category of $\infty$-prestacks equipped with the (global) projective model structure, then the model category constructed via left Bousfield localization at the class $\mathcal{S}$ of all hypercovers \cite[section 2]{DHI}, denoted $\sPre(\Bihol)_{/{\mc S}}$, will be referred to as the \emph{local projective model structure}.
\end{definition}

\begin{definition}
An \emph{$\infty$-stack} is a fibrant $\infty$-prestack in the local projective model structure $\sPre(\Bihol)_{/{\mc S}}$.
\end{definition}

While the context of this paper conveniently avoids the need to explicitly use, and so the need to provide a description of, hypercovers, the \v{C}ech nerve of an covering family is sufficient for all purposes here. Traditionally, when $\mathscr{C}$ is the site of open covers $\{U_i \subset X\}_{i \in I}$ for a fixed topological space $X$, the fibered products $U_i \times_X U_j$ are exactly equal to the intersection $U_i \cap U_j$. From here the \v{C}ech nerve as a simplicial presheaf is defined by taking coproducts of representable presheaves induced by the higher intersections $U_{i_0} \cap \cdots \cap U_{i_p}$. However for the covering families $\mc U = \{U_i \xrightarrow{f_i} W\}_{i \in I}$ in $\Bihol$, the pullback of a pair of morphisms has two canonical presentations:
\begin{equation}
\begin{tikzcd}
f_j^{-1}\left(f_i(U_i) \cap f_j(U_j)\right)\cong f_i^{-1}\left(f_i(U_i) \cap f_j(U_j)\right) \arrow[r, dotted] \arrow[d, dotted]& U_i \arrow[d, "f_i"] \\
U_j \arrow[r, "f_j"]& W.
\end{tikzcd}
\end{equation}

To define the \v{C}ech nerve below, a consistent choice is needed given any family of morphisms $\{U_{i} \xrightarrow{f_{i}} U\}_{i \in I}$ in $\Bihol$:
\begin{equation}\label{EQ Cech pullback choice}
U_{i_0, \ldots, i_p}:= f_{i_0}^{-1} \left( f_{i_0}\left(U_{i_0}\right) \cap \cdots \cap f_{i_p}\left(U_{i_p}\right)  \right) \cong U_{i_0} \times_U \cdots \times_U U_{i_p}.
\end{equation}

There are two uses for defining a \v{C}ech nerve in this paper: first is to check that certain $\infty$-prestacks are stacks and the second is for cofibrant replacement of certain $\infty$-prestacks. The former only requires the \v{C}ech nerve to be defined for a covering family in $\text{Cov}(\Bihol)$ while the latter requires a slightly more general construction for the analogous family of morphisms covering a general $n$-dimensional complex manifold. 
\begin{definition}\label{DEF covering family}
Let $\mathcal{U} =\{U_i \xrightarrow{f_i} U\}_{i \in I}$ be a covering family in $\Bihol$. The \emph{\v{C}ech nerve of a covering family $\mathcal{U}$} is the $\infty$-prestack $\Cech\mathcal{U}\in \sPre(\Bihol)$ whose set of $p$-simplices on a test object $W$ in $\Bihol$ is defined to be the coproduct of representable presheaves of sets:
\begin{equation}
\Cech\mathcal{U}(W)_p:= \coprod\limits_{i_0, \ldots, i_p} \yo U_{i_0, \ldots, i_p}(W) 
\end{equation}
where $U_{i_0, \ldots, i_p}$ follows \eqref{EQ Cech pullback choice} and $\yo U_{i_0, \ldots, i_p}(W):= \Bihol(W,U_{i_0, \ldots, i_p})$ is the Yoneda embedding $\left(\Bihol \right)^{op} \xrightarrow{\yo} \sPre(\Bihol)$. 
\end{definition}
Note that the above coproduct can be thought of as being taken over all non-empty fibered products since the Yoneda embedding applied to the empty manifold will be the empty set, whose coproduct-component then adds nothing to the $p$-simplices of the \v{C}ech nerve.

\begin{proposition}\label{PROP stack cond}
An $\infty$-prestack $\mc F$ on $\Bihol$, which objectwise takes values in Kan complexes of bounded homotopy type, is an $\infty$-stack if for each covering family $\mc U$ of $U$ the map 
\begin{equation}\label{EQ: descent TOT}
\mathcal{F}(U)  \to \Tot\left(\mathcal{F}_{\bullet} \left(  \Cech\mathcal{U}^{\bullet} \right) \right)
\end{equation}
is a weak equivalence, where the codomain is the totalization \cite[definition 18.6.3]{Hi} of the cosimplicial simplicial set $\mathcal{F}_{n} \left(  \Cech\mathcal{U}^{p} \right):= \prod\limits_{i_0, \ldots, i_p} \mathcal{F}_{n}\left( U_{i_0, \ldots, i_p}\right)$.
\begin{proof}
By \cite[theorem 1.3]{DHI} the fibrant objects of $\sPre(\Bihol)_{/{\mc S}}$ are those which are objectwise fibrant as simplicial sets and satisfy descent with respect to all hypercovers. But then by \cite[corollary 8.6]{DHI} since $\mc F$ has bounded homotopy type it will satisfy descent with respect to all hypercovers if it satisfies descent with respect to all \v{C}ech covers. In other words, for any map $\Cech\mathcal{U} \to U$ the induced map on hom-spaces, 
\begin{equation}\label{EQ: descent}
\underline{\sPre}({\Bihol})(\yo U, \mathcal{F} ) \to  \underline{\sPre}({\Bihol})(\Cech \mathcal{U} , \mathcal{F})
\end{equation}
needs to be a weak equivalence. As usual, by a consequence of the \v{C}ech nerve being Reedy cofibrant \cite[lemma C.5]{GMTZ2}, the right hand side of \eqref{EQ: descent} is computed via totalization and the left hand side always reduces to $\mathcal{F}(U)$  by the simplicial Yoneda lemma (see  \cite[section VIII.1, lemma 1.2]{GJ} or \cite[eq 2.31]{K}) so that one only needs to check \eqref{EQ: descent TOT} as stated.
\end{proof}
\end{proposition}

\section{An $\infty$-stack of charts}
In this section, we define the $\infty$-stack $\Chart:(\Bihol)^{op}\to \sSet$. Given an object $W$ in ${\Bihol}$, consider the simplicial set $\Chart (W)$ whose $k$-simplices, consist of commutative diagrams generated by the nerve of the diagram
 \begin{equation}\label{EQU:phi_{j,j+1}-diagram}
\begin{tikzcd}
&&W\arrow[swap]{dll}{\rho_0}\arrow{dl}{\rho_1}\arrow[swap]{dr}{\rho_{k-1}}\arrow{drr}{\rho_k} &&\\
V_0 & V_1\arrow{l}{\phi_{0,1}}  & \dots \arrow{l}{\phi_{1,2}} & V_{k-1} \arrow{l}{\phi_{k-2,k-1}} & \arrow{l}{\phi_{k-1,k}} V_k
\end{tikzcd}
\end{equation}
where each $V_i$ is an open subset of $\mathbb{C}^n$, each vertical arrow $\rho_i$ is a biholomorphism, and for each $1 \le i \le j \le \ell$ there is a unique biholomorphism $\phi_{i,j}=\rho_i\circ \rho_j^{-1}: V_j \to V_i$ making the diagram commute.  Such diagrams pull back along any map $W' \xrightarrow{\phi} W$ in ${\Bihol}$ by pullback: 
\begin{equation}
\begin{tikzcd}[row sep = small, column sep=small]
 & & W' \arrow[hook]{rr}{\phi}&  & W  \\
& &  \Chart(W')& & \Chart(W)  \arrow{ll}[swap]\Chart(\phi)\\
&W' \arrow[""name=x]{ddl}[swap]{\rho_0 \circ \phi}   \arrow{ddr}{\rho_k\circ \phi} & & &  && W \arrow[""name=y]{ddl}[swap]{\rho_0}  \arrow{ddr}{\rho_k} \\
&&  \phantom{A} & &   \phantom{B}\arrow[mapsto]{ll}&\\
(\rho_0 \circ \phi)(W')&\arrow{l}{\phi_{0,1}}\cdots &\arrow{l}{\phi_{k-1,k}}(\rho_k \circ \phi)(W') &  && V_0   &\arrow{l}{\phi_{0,1}}\cdots &V_k\arrow{l}{\phi_{k-1,k}}\label{EQ: Charts pullback}
\end{tikzcd}
\end{equation}   
so we obtain an $\infty$-prestack $\Chart$ on $\Bihol$. Note in particular that the only effect that $\Chart(\phi: W' \to W)$ has on each $\phi_{i,j}$ is to restrict it to the appropriate domain $(\rho_j \circ \phi)(W')$. A more formal definition is now offered.  

\begin{definition}\label{DEF pre Chart_0}
The \emph{$\infty$-prestack of charts}, $\Chart: \left({\Bihol}\right)^{op} \to \sSet$ is defined:
\begin{itemize}
\item  for each object $W \in \Bihol$ to be the $0$-coskeletal simplicial set $\Chart(W)$ whose vertex-set $\Chart(W)_0$ is equal to the set of all biholomorphisms to open subsets of $\mathbb{C}^n$, with source $W$, i.e., $\Chart(W)_0:=\left\{ W \xrightarrow{\cong} V \subset \mathbb{C}^n \right\}$, 
\item for each morphism $\phi: W' \to W$ (i.e., an open embedding which is a biholomorphism onto its image), the associate map between vertex sets is the pullback along $\phi$, $\Chart(\phi)_0(\rho):= \rho \circ \phi$, with the codomain appropriately restricted to result in a surjection (see \eqref{EQ: Charts pullback}).
\end{itemize}
\end{definition}
Using \eqref{EQU:phi_{j,j+1}-diagram} as a guide, a $k$-simplex in $\Chart(W)_k$ is such a diagram where the vertical maps, i.e., the vertices, uniquely determine the rest of the simplex.

\begin{lemma}\label{LEM: cosk tot}
If $X^{\bullet}_{\bullet}$ is a cosimplical simplicial set which is (cosimplicial) degree-wise a $0$-coskeletal simplicial set then so is its totalization.
\begin{proof}
As $X^{m}_{\bullet}$ is assumed to be $0$-coskeletal for each $m \ge 0$ then any collection of $(k+1)$-many vertices $\nu = \{v_0, \ldots, v_k\}$ uniquely determines a $k$-simplex $\widetilde{\nu} \in X^{m}_{k}$. So let $\alpha: \partial \Delta^m \to \Tot\left( X^{\bullet}_{\bullet}\right)$ be a map. Recall that an $m$-simplex of $\Tot\left( X^{\bullet}_{\bullet}\right)$ would consist of, for each $p\ge 0$, a map $\widetilde{\alpha^p}: \Delta^m \times \Delta^p \to X^p_{\bullet}$ satisfying the usual coherence conditions. But $\alpha$ in particular begins with $(m+1)$-many vertices $\nu^0 = \{v_0^0, \ldots, v^0_m\}$ in $X^0_{0}$. Since $X^0$ is coskeletal, these vertices uniquely determine the required map $\widetilde{\alpha^0}: \Delta^m \times \Delta^0\to X^0_{\bullet}$. Propogating these vertices to $p(m+1)$ many vertices in $X^{p-1}_{\bullet}$, using the cosimplicial maps, provides the required vertices of $\alpha$ at cosimplicial level $p$ due to the coherence conditions of totalization. The fact that each $X^p$ is coskeletal then means that these vertices uniquely determine $\widetilde{\alpha}$ in a way that agrees with ${\alpha}$ at each cosimplicial degree.
\end{proof}
\end{lemma}

\begin{proposition}\label{PROP chart stack}
${\Chart}$ is an $\infty$-stack on ${\Bihol}$.
\begin{proof}
$\Chart$ is object-wise fibrant, i.e., ${\Chart}(W)$ is a Kan complex for each $W$ since every $0$-coskeletal space is contractible. Moreoever, by lemma \ref{LEM: cosk tot} $\Tot\left({\Chart}\left(  \mathcal{U} \right) \right)$ is $0$-coskeletal and so also contractible. Since the map \eqref{EQ: descent TOT} is now a map between contractible spaces in this case then it is a weak equivaleance.
\end{proof}
\end{proposition}

\section{Two $\infty$-prestacks of forms}

We now define two $\infty$-prestacks $\extAnti$ and $\holFormP{k}$ and the map $\Hartmap: \extAnti \to \holFormP{n}$ of $\infty$-prestacks, which comes from the Hartogs extension theorem.
For a complex manifold $M$ there is the classic bi-graded $(\del, \delbar)$ complex $A^{(\bullet, \bullet)}(M) = \bigoplus\limits_{p,q \ge 0}A^{(p,q)}(M)$ of forms on $M$ of (holomorphic, antiholomorphic) bi-degree $(p,q)$. The graded subspace of holomorphic forms is denoted $\holForm^{\bullet}(M)\subset A^{(\bullet, 0)}(M)$.  On any open subset of a product manifold $N \subset M\times M'$ we can further decompose the complex into components $A^{(p,q),(p',q')}(N)$ where $p$ and $p'$ are the holomorphic grading on $M$ and $M'$ respectively whereas $q$ and $q'$ are the anti-holomorphic grading on $M$ and $M'$ respectively. Of particular interest for this paper is the case where $M'=M$, we further consider germs of forms near the diagonal $A_{\Delta}^{(p,q),(p',q')}(M):= \lim\limits_{N\supset \Delta}A^{(p,q),(p',q')}(N - \Delta)$ where the limit runs over open neighborhoods $N$ of the diagonal $\Delta \subset N \subset M \times M$ under inclusions. The motivation for the following construction is to provide a receptacle for Toledo and Tong's \v{C}ech parametrix $\omega^0, \ldots, \omega^{n-1}$, along with corresponding holomorphic cochain $\omega^n$ to form a single cocycle.
\begin{definition}\label{DEF extAnti}
Given an $n$-dimensional complex manifold $M$, define an augmented \emph{$\delbar$-complex of forms} written 
\begin{equation}
\extAnti^{\bullet}(M)= \left( 0 \to \extAnti^{-n}(M)\hookrightarrow \extAnti^{1-n} (M)\xrightarrow{\delbar}  \extAnti^{2-n}(M) \xrightarrow{\delbar} \dots   \xrightarrow{\delbar} \extAnti^{0} (M) \to 0 \right)
\end{equation} to be the graded $\holForm^0(M)$-module, 
\begin{equation}\label{EQ extAnti def}
\extAnti^p(M):= \begin{cases} \bigoplus\limits_{q+q'=n-1} Z_{\Delta}^{(0,q),(n,q')}(M \times M - \Delta)&\quad \text{ if } p=0 \\[1em]
\bigoplus\limits_{q+q'=n-1+p} A_{\Delta}^{(0,q),(n,q')}(M \times M - \Delta) &\quad \text{ if } p=1-n, \ldots,- 1 \\[1.5em]  \lim\limits_{N \supset \Delta }\holForm^{(0,n)}(N - \Delta)  &\quad \text{ if } p=-n \\ 0 &\quad \text{otherwise}
\end{cases}
\end{equation}
where $Z_{\Delta}^{(0,q),(n,q')}(M)\subset A_{\Delta}^{(0,q),(n,q')}(M)$ is the subspace of germs of $\delbar$-closed forms, and the differential is $\delbar$ in all degrees except for $\extAnti^{-n}(M)\hookrightarrow \extAnti^{1-n} (M)$ which is the inclusion of those germs of holomorphic $n$-forms near the diagonal that live in $A_{\Delta}^{(0,0),(n,0)}(M)$.
\end{definition}

We want to apply the Dold-Kan functor to $\extAnti^{\bullet}(M)$. To be precise and to fix our notation, we next give an  explicit description of the Dold-Kan functor.
\begin{notation}\label{DEF:DK(Cbullet)}
Let $ \left(\mc{C}^{\bullet}, d\right)$ be a non-positively graded cochain complex. Applying the Dold-Kan functor gives a simplicial set $\DK(\mc{C}^{\bullet})=Hom(N(\Z\Delta^\bu ),\mc{C}^{\bullet})$, see \cite[Appendix B]{GMTZ1} for the notation used here, which is explicitly given as follows:

If we denote the $\ell$-cells of the standard $n$-simplex by $e_{i_0,\dots, i_\ell}$ for $0\leq i_0<\dots<i_\ell\leq n$, then an $n$-simplex in $\DK(\mc{C}^{\bullet})_n$ is a labeling $\{c_{i_0,\dots,i_\ell} \in \mc{C}^{-\ell} \}_{i_0,\dots,i_\ell}$ of all cells $e_{i_0,\dots, i_\ell}$ by elements $c_{i_0,\dots, i_\ell}$, where we require that the alternating sum of the boundary of any $(k+1)$-cell is in the image of the  differential $d$ of $\mc{C}$ applied to the interior:
 \begin{equation}\label{EQ DK or dgnerve}
 \sum_{j=0}^{k+1} (-1)^j c_{i_0,\dots, \widehat{i_j},\dots,i_{k+1}}=d\left(c_{i_0,\dots,i_{k+1}} \right).
 \end{equation} 

The simplicial face map $d_j:\DK(\mc{C}^{\bullet})_n\to \DK(\mc{C}^{\bullet})_{n-1}$ is given by labeling $d_j(\{c_{I}\})$ by the labelings of the $j^\text{th}$ face of $\{c_{I}\}_I$, i.e., we set $d_j(\{c_{I}\})_{i_0,\dots, i_\ell}=c_{\delta^j(i_0),\dots, \delta^j(i_\ell)}\in C^{-\ell}$, where $\delta^j:\{0,\dots,n-1\}\to \{0,\dots, n\}$ is the map that skips $j$. 

Similarly, the simplicial degeneracy map $s_j:\DK(\mc{C}^{\bullet})_n\to \DK(\mc{C}^{\bullet})_{n+1}$ is given by setting $s_j(\{c_{I}\})_{i_0,\dots, i_\ell}=c_{\sigma^j(i_0),\dots, \sigma^j(i_\ell)}$, where $\sigma^j:\{0,\dots,n+1\}\to \{0,\dots, n\}$ repeats $j$. Here, labels with repeated indices $c_{\dots,j,j,\dots}$ are set to be $0$, and thus, we only obtain non-zero values for the $j^\text{th}$ and $(j+1)^\text{st}$ face of $s_j(\{v_{I}\})$, where one of the double indexing by $j$ is removed. {For example, the element $s_j(\{c_{I}\})_{i_0,\dots, i_\ell}\in C^{-\ell}$, evaluated in the case where we take the indices $(i_0,\dots, i_\ell)=(0,\dots,\widehat{j},\dots,n+1)$, becomes 
\begin{equation}\label{EQU:degen-sj(cI)}
s_j(\{c_{I}\})_{0,\dots,j-1,\widehat{j},j+1,\dots, n+1}=c_{\sigma^j(0),\dots, \sigma^j(j-1),\widehat{\sigma^j(j)},\sigma^j(j+1),\dots, \sigma^j(n+1)}\\=c_{0,\dots,j-1,\widehat{j},j,\dots, n}=c_{0,\dots,n},
\end{equation}
which does not have to be zero.}
\end{notation}

\begin{example}[Dold-Kan of a vector space in degree $k$]\label{DEF:DK(V)}
Let $V$ be a vector space and $k\geq 0$. Let $V[-k]$ denote shifting $V$ down by $k$, and, by abuse of notation, write $V[-k]$ as well for the cochain complex with $0$ in all other degrees $\neq -k$. 
An $n$-simplex in $\DK(V[-k])_n$ is, again by notation \ref{DEF:DK(Cbullet)}, a labeling $\{v_{i_0,\dots,i_\ell}\}_{i_0,\dots,i_\ell}$ of all cells $e_{i_0,\dots, i_\ell}$ but where $v_{i_0,\dots,i_\ell}=0$ for $\ell\neq k$, and $v_{i_0,\dots, i_k}\in V$ are such that the alternating sum of the boundary of any $(k+1)$-cell vanishes, $\sum_{j=0}^{k+1} (-1)^j v_{i_0,\dots, \widehat{i_j},\dots,i_{k+1}}=0$. 
\end{example}

Now applying the Dold-Kan functor in two different cases will help construct the two $\infty$-prestacks of differential forms in this paper, $\extAnti$ (definition \ref{DEF extAnti pre}) and $\holFormP{n}$ (definition \ref{DEF:Okhol}):

\begin{definition}\label{DEF extAnti pre}
Define the $\infty$-prestack $\extAnti: \left(\Bihol\right)^{op}\to \sSet$ by setting $\extAnti(W):= DK  \left( {\extAnti^{\bullet}}(W)\right)$. For a morphism $\phi:W'\to W$ in $\Bihol$, we set $\extAnti(\phi)_n:\extAnti(W)_n\to \extAnti(W')_n$ to be the pullback of labelings $c_{i_0,\dots, i_\ell}\in\extAnti^{-\ell}(W)_n$, i.e., $\extAnti(\phi)_n(\{c_I\})_{i_0,\dots, i_\ell}:=\phi^*(c_{i_0,\dots, i_\ell})$.
\end{definition}
\begin{definition}\label{DEF:Okhol}
For a fixed $k>0$, define the simplicial presheaf $\holFormP{k}:\left(\Bihol\right)^{op}\to \sSet$ by setting $\holFormP{k}(W):=DK(\holForm^k(W)[-k])$. For a morphism $\phi:W'\to W$ in $\Bihol$, we set $\holFormP{k}(\phi)_n:\holFormP{k}(W)_n\to \holFormP{k}(W')_n$ to be the pullback of labelings $v_{i_0,\dots, i_\ell}\in \holForm^k(W)_n$, i.e., $\holFormP{k}(\phi)_n(\{v_I\})_{i_0,\dots, i_\ell}:=\phi^*(v_{i_0,\dots, i_\ell})$.
\end{definition}

\begin{proposition}\label{PROP holForm stack}
The $\infty$-prestack $\holFormP{k}$ on ${\Bihol}$ is an $\infty$-stack. 
\end{proposition}
In order to prove that holomorphic $k$-forms induce an $\infty$-prestack the following lemma, relating the Dold-Kan functor and the totalization functors, is useful. The statement requires a ``smart truncation'' functor, $\tau_{\le 0}$  which for a $\mathbb{Z}$-graded complex $C^{\bullet}$ assigns \cite[Truncations 1.2.7]{We} ,
\begin{equation}\label{EQ smart trunc}
\left(\tau_{\le 0}C^{\bullet} \right)^p= \begin{cases}C^p &\text{ if } p< 0\\ Z^0 & \text{ if }p=0  \\ 0 &\text{ if } p > 0 \end{cases}.
\end{equation}
In all of the applications in this paper, (see for example \eqref{EQ extAnti def}) the complexes are already ``pre-truncated''. 
\begin{lemma}\label{LEM Emilios Tot DK}
If $\mc A$ is a presheaf of chain complexes and $\mc U$ is a cover of $W$ then 
\begin{equation}\label{EQ Emilios Tot DK}
\Tot \left((\DK \circ \tau)\mc A \left(  \Cech\mathcal{U}^{\bullet} \right)\right)
{=} \DK   \left( \left( \tau \circ \tot\right)   \mc A  \left(  \Cech\mathcal{U}^{\bullet}\right) \right)
\end{equation}
\begin{proof}
See \cite[eqn. (113)]{Mi}.
\end{proof}
\end{lemma}

\begin{proof}[Proof of proposition \ref{PROP holForm stack}]
Applying proposition \ref{PROP stack cond} requires verifying three conditions: (i) $\holFormP{k}$ is object-wise Kan, (ii) $\holFormP{k}$ is of bounded homotopy type, and (iii) the map \eqref{EQ: descent TOT} is a weak equivalence when $\mc F = \holFormP{k}$.

The first is satisfied since simplicial abelian groups are Kan complexes \cite[page 12]{GJ} and $\holFormP{k}$ is formed object-wise by taking Dold-Kan, landing in simplicial abelian groups. The third condition amounts to checking that the map \eqref{EQ: descent TOT} in the case of $\mc F = \holFormP{k}$: 
\begin{align*}
\holFormP{k}(W)=\DK  \left( \holFormP{k}[-k]\right) (W) \to \Tot\left(\holFormP{k} \left(  \Cech\mathcal{U}^{\bullet} \right) \right)=  &\Tot \left(\DK  \left( \holFormP{k}[-k]\right) \left(  \Cech\mathcal{U}^{\bullet} \right)\right)\\
\stackrel{\text{\eqref{EQ Emilios Tot DK} }}{=}&  \DK   \left( \left( \tau \circ \tot\right)  \left(  \holFormP{k}[-k]^{\bullet} \left(  \Cech\mathcal{U}^{\bullet}\right) \right)\right)
\end{align*}
is a weak equivalence.  Conveniently, both sides are given as the Dold-Kan of some complex, and so, since the homotopy groups of Dold-Kan are the homology of the underlying complexes, it remains to prove that the associated map:
 \[ H^n \left(  \holFormP{k}[-k](W) \right) \to H^n  \left( \left( \tau \circ \tot\right)  \left(  \holFormP{k}[-k]^{\bullet} \left(  \Cech\mathcal{U}^{\bullet}\right) \right)\right) \]
 is an isomorphism. The left hand side is only non-trivial in degree $-k$ where it is equal to $\holForm^{k} (W)$, which covers the second condition that $\holFormP{k}$ has bounded homotopy type. Turning to the right side, the bi-complex $\holFormP{k}[-k]^{\bullet} \left(  \Cech\mathcal{U}^{\bullet}\right)$ is trivial in every ``row'' and so 
 \[H^n  \left( \left( \tau \circ \tot\right)  \left(  \holFormP{k}[-k]^{\bullet} \left(  \Cech\mathcal{U}^{\bullet}\right) \right)\right) = H^n \left(   \holFormP{k} \left( \Cech\mathcal{U}^{\bullet}\right)[-k] \right)=H^{n+k} \left(   \holFormP{k} \left( \Cech\mathcal{U}^{\bullet}\right) \right). \]
As a presheaf of $\mathcal{O}_{(-)}$-modules, $\holForm^{k}$ is coherent analytic (as it is locally-free) implying that, by Cartan's Theorem\footnote{Note this is the sole part of the paper where the Stein property is used.} B, the sheaf cohomology $H^p(\holForm^{k})$ is trivial in positive degrees. Hence $H^{n+k} \left(   \holFormP{k} \left( \Cech\mathcal{U}^{\bullet}\right) \right)$ is trivial except for $n=-k$ in which case it computes as $\holForm^{k} (W)$.
\end{proof}
We now use the Hartogs extension theorem \cite[theorem (3.28)]{De} to define a map of infinity prestacks, $\Hartmap$.
\begin{proposition}\label{PROP Hartmap is a map}
The following assignment forms a map of $\infty$-prestacks $\Hartmap: \extAnti \to \holFormP{n}$: For each $n$-simplex labeled on the top cell by $\alpha^n\in \extAnti^{-n}(W)=  \lim\limits_{N \supset \Delta}\holForm^{(0,n)}(N-\Delta)$, assign the holomorphic form $\Delta^*{\overline{\alpha^n}} \in \holForm^{n}(W)$ which is the holomorphic extension $\overline{\alpha^n}$  of any representative of the germ to $N$, which is then restricted along the diagonal $W\cong \Delta \hookrightarrow W \times W$. For $k \ne n$ assign the trivial $k$-simplex.
\begin{proof}The assignment $\Hartmap$ is actually Dold-Kan applied to an assignment $\widetilde{\Hartmap}: \extAnti^{\bullet}(-) \to \holFormP{n}(-)^{\bullet}$ of presheaves of cochain complexes:
\begin{equation}
\widetilde{\Hartmap}\left(\omega \in  \extAnti^{p}(W) \right) := \begin{cases} \Delta^*\left(\overline{\omega} \right) \in \holForm^{n}(W) &\text{ if } p=-n\\ 0 &\text{ otherwise,} \end{cases}\label{EQ Hartog Delta}
\end{equation}
where $\overline{\omega}$ is the holomorphic extension of a representative of $\omega$ on $N - \Delta $ to all of $N$. Thus, it suffices to check $\widetilde{\Hartmap}$ is a map of presheaves. For each $W$, $\widetilde{\Hartmap}(W)$ is trivially a chain map due to the codomain being concentrated in degree $-n$. It is a map of presheaves since pullback to the diagonal commutes with pullback along a biholomorphism (and extension is unique). 
\end{proof}
\end{proposition}

\begin{definition}\label{DEF hartogs}
The map of $\infty$-prestacks $\Hartmap: \extAnti \to \holFormP{n}$ will be referred to as the \emph{Hartogs map}.
\end{definition}

\section{Toledo and Tong's \v{C}ech parametrix as a map of $\infty$-prestacks}\label{SECTION TT parametrix}
The universal \v{C}ech parametrix constructed by Toledo and Tong in \cite{TT76}, can be reinterpreted as a map of $\infty$-prestacks $\BMmap: \Chart \to \extAnti$, as described in this section. Begin with the ($\delbar$-closed) Bochner-Martinelli kernel, $\omega ^0\in \extAnti^{0}(\mathbb{C}^n)$, defined by
\begin{equation}\label{DEF Bochner-Martinelli kernel}
\omega^{0}(z, \xi) =b_n \frac{\sum_{i=1}^n (-1)^i \left( \overline{\xi^i} -  \overline{z^i}\right) \cdot \left( d\overline{\xi^1} -  d\overline{z^1}\right) \cdots \widehat{ \left( d\overline{\xi^i} -  d\overline{z^i}\right)}  \cdots \left( d\overline{\xi^n} -  d\overline{z^n}\right)}{\abs{\xi - z}^{2n}} d\xi^1 \dots d\xi^n,
\end{equation}
where $b_n: =  -(-1)^{n(n-1)}\frac{(n-1)!}{2\pi i}$. For each vertex $\rho: W \to V\subset \C^n$ in $\Chart(W)_0$, restricting $\omega^0$ to $V$ and then pulling back along $\rho$ provides the object-wise map on vertices:
\begin{align}
 \Chart(W)_0 &\xrightarrow{\BMmap(W)_0} \extAnti(W)_0\\
 \left( W \xrightarrow{\rho} V\right) &\xmapsto{\phantom{\BMmap(W)_0}} (\rho \times \rho)^*\omega^0. \nonumber
\end{align}
For each $k$-simplex $\rho_{\bullet}= (\rho_0:W\to V_0,\dots, \rho_k:W\to V_k) \in \Chart(W)_k$, as in \eqref{EQU:phi_{j,j+1}-diagram}, $\BMmap(W)_k$ needs to assign a $k$-simplex in $\extAnti(W)_k = DK  \left( \extAnti^{\bullet}(W)\right)_k$, as in definition \ref{DEF extAnti pre}. Recall from notation \ref{DEF:DK(Cbullet)} that such a $k$-simplex amounts to a labeling of each $\ell$-cell $e_{i_0,\dots, i_\ell}$ given by $0 \le i_0 < \cdots < i_{\ell} \le k$, by an element $\left(\BMmap(W)_k(\rho_{\bullet})\right)_{i_0, \ldots, i_{\ell}}$ in $\extAnti^{-\ell}(W)$ satisfying
 \begin{equation}\label{EQ BM sat DK}
d\left(\left(\BMmap(W)_k(\rho_{\bullet})\right)_{i_0, \ldots, i_{\ell}} \right)=  \sum_{j=0}^{\ell} (-1)^j \left(\BMmap(W)_k(\rho_{\bullet})\right)_{i_0, \ldots,\widehat{i_{j}}, \ldots,  i_{\ell}}.
 \end{equation} 
 where $d$ is the differential from definition \ref{DEF extAnti}.  
 
In order to define $\BMmap$ on higher simplices we will make use of the universal \v{C}ech parametrix that was constructed in \cite{TT76}. Recall from \cite[definition 4.3 and the paragraph following it]{TT76} that the universal \v{C}ech parametrix provides a sequence $(\omega^0, \omega^1, \ldots, \omega^{n-1})$ with $\om^q(\phi_1,\dots, \phi_q)\in \extAnti^{-q}(dom(\phi_q))$, where $(\phi_1, \ldots, \phi_q)$ is a string of $q$ composable morphisms (for $q<n$), which starts with the Bochner-Martinelli kernel $\om^0$, and which satisfies the condition $\delbar \omega^{q} = \delta \omega^{q-1}$ from \cite[(4.4)]{TT76}, i.e., 
\begin{multline}\label{EQ TT cocycle}
\delbar \omega^{q}\left( \phi_1, \ldots, \phi_{q}\right) = \delta(\omega^{q-1})\left( \phi_1, \ldots, \phi_{q}\right)
\\= \omega^{q-1}\left( \phi_2, \ldots, \phi_{q}\right) + \sum\limits_{j=1}^{q-1} (-1)^j \omega^{q-1}\left( \phi_1, \ldots,\phi_{j}\circ \phi_{j+1}, \ldots,  \phi_{q}\right)+\left(\phi_{q}\times \phi_{q}\right)^*\omega^{q-1} \left( \phi_1, \ldots, \phi_{q-1}\right)
\end{multline}
for all $q=1,\dots,n-1$. Moreover, we recall that we also have the associated holomorphic cochain $\omega^n(\phi_1,\dots, \phi_n)\in  \extAnti^{-n}(dom(\phi_n))$ defined by $\omega^n:= \delta \omega^{n-1}$. 
The existence of such a universal parametrix is provided by \cite[theorem 10.9]{TT76}.

\begin{definition}\label{DEF:BM-on-simplicies}
Let $W$ be an object in $\Bihol$ and let $\rho_{\bullet}= (\rho_0:W\to V_0,\dots, \rho_k:W\to V_k) \in \Chart(W)_k$ be a $k$-simplex. We use the notation $\phi_{i,j}:= \rho_i \circ \rho_j^{-1}$ as before, and for indices $j_0,\dots, j_\ell\in\{0,\dots, k\}$, we set
$$\omega_{\rho_{j_0},\dots,\rho_{j_\ell}}:= \left( \rho_{j_\ell} \times \rho_{j_\ell}\right)^*\left(\omega^{\ell}(\phi_{j_0,j_1}, \ldots, \phi_{j_{\ell -1},j_{\ell}}) \right),$$ which is the \v{C}ech parametrix applied to transition function coming from the charts $(\rho_{j_0},\dots,\rho_{j_\ell})$ and which is then pulled back to $W$. We then define $\BMmap(W):\Chart(W)\to \extAnti(W)$ via a skew-symmetrized version of these $\omega$'s. More precisely, on $k$-simplicies $\BMmap(W)_k:\Chart(W)_k\to \extAnti(W)_k$,  $\rho_\bu \mapsto\BMmap(W)_\ell(\rho_\bu)$, we label $\BMmap(W)_k(\rho_\bu)_{i_0,\dots, i_\ell}$ for each $\ell$-cell given by some indices $0 \le i_0 < \cdots < i_{\ell} \le k$ as follows:
\begin{itemize}
\item if $\ell\leq n$, we define 
\begin{align}\label{EQ define BM cells}
\left(\BMmap(W)_k(\rho_{\bullet})\right)_{i_0, \ldots, i_{\ell}}:=&\frac{1}{(\ell +1)! }\sum_{\sigma_\in S_{\ell+1}} \text{sgn}(\sigma)\cdot \omega_{\rho_{i_{\sigma(0)}},\dots,\rho_{i_\sigma(\ell)}}
\end{align}
\item if $\ell>n$, then we set $(\BMmap(W)_k(\rho_{\bullet}))_{i_0, \ldots, i_{\ell}}:=0$
\end{itemize}

\end{definition}
We claim that this definition gives a map of $\infty$-prestacks.
\begin{proposition}\label{PROP TTmap is a map}
The above definition \ref{DEF:BM-on-simplicies} gives a map of $\infty$-prestacks $\BMmap: \Chart \to \extAnti$.
\begin{proof}
We show that $\Bihol$ is a natural transformation in three steps. First, we show that the above assignment \eqref{EQ define BM cells} gives a well-defined $k$-simplex, second, the assignment $\BMmap(W)$ respects face and degeneracy maps, and third, $\BMmap$ gives the correct commutative diagram for a natural transformation.
\begin{itemize}
\item[\bf Step 1:\hspace{-2.5mm}] \hspace{1.3mm}
First, we check that for a fixed object $W$ in $\Bihol$, the forms defined in \eqref{EQ define BM cells} give well-defined $\ell$-cells, i.e., that they satisfy the ``Dold-Kan condition'' from \eqref{EQ BM sat DK}. To check \eqref{EQ BM sat DK} for $\ell \ne n$ on the un-symmetrized data, we use the \v{C}ech parametrix condition \eqref{EQ TT cocycle} as follows:
\begin{align}
&\delbar\left( \left(\rho_{{i_{\ell}}}\times \rho_{{i_{\ell}}}\right)^*\left(\omega^{\ell}(\phi_{i_0,i_1}, \ldots, \phi_{i_{\ell -1},i_{\ell}}) \right)\right)\label{LHS BM DK}\\
= & \left(\rho_{{i_{\ell}}}\times \rho_{{i_{\ell}}}\right)^* \left(\delbar \left(\omega^{\ell}(\phi_{i_0,i_1}, \ldots, \phi_{i_{\ell -1},i_{\ell}}) \right)\right)\nonumber\\
\substack{\eqref{EQ TT cocycle}\\=}&\left(\rho_{{i_{\ell}}}\times \rho_{{i_{\ell}}}\right)^* \left( \omega^{\ell-1}(\phi_{i_1,i_2}, \ldots, \phi_{i_{\ell -1},i_{\ell}}) \right)\nonumber
\\
&+\left( \rho_{{i_{\ell}}} \times \rho_{{i_{\ell}}}\right)^* \left( \sum\limits_{j=1}^{\ell-1} (-1)^j \omega^{\ell-1}\left( \phi_{i_0,i_1}, \ldots,\phi_{i_{j-1},i_j}\circ \phi_{i_j,i_{j+1}}, \ldots,  \phi_{i_{\ell -1},i_{\ell}}\right)\right)\nonumber\\
&+\left( \rho_{{i_{\ell}}}\times \rho_{{i_{\ell}}}\right)^* \left( \left( \phi_{i_{\ell -1},i_{\ell}}\times \phi_{i_{\ell -1},i_{\ell}}\right)^* \left(\omega^{\ell-1}(\phi_{i_0,i_1}, \ldots, \phi_{i_{\ell -2},i_{\ell-1}}) \right)\right)\nonumber\\
=&\left( \rho_{{i_{\ell}}}\times \rho_{{i_{\ell}}}\right)^* \left(\omega^{\ell-1}(\phi_{i_1,i_2}, \ldots, \phi_{i_{\ell -1},i_{\ell}}) \right)
)\label{RHS BM DK}\\
&+ \sum\limits_{j=1}^{\ell-1} (-1)^j \left(\rho_{{i_{\ell}}}\times \rho_{{i_{\ell}}} \right)^* \omega^{\ell-1}\left( \phi_{i_0,i_1}, \ldots,\phi_{i_{j-1},i_{j+1}}, \ldots,  \phi_{i_{\ell -1},i_{\ell}}\right)\nonumber\\
&+\left(\rho_{\ell-1}\times \rho_{\ell-1} \right)^*  \left(\omega^{\ell-1}(\phi_{i_0,i_1}, \ldots, \phi_{i_{\ell -2},i_{\ell-1}}) \right).\nonumber
\end{align}
Since skew-symmetrization commutes with $\delta$ as well as with $\delbar$ (see \cite[pg 299]{TT76}), the first claim is verified for $\ell < n$ by observing that the left and right hand sides of \eqref{EQ BM sat DK} are the skew-symmetrizations of \eqref{LHS BM DK} and \eqref{RHS BM DK}, respectively.

For $\ell =n$ the holomorphic cochain is defined via $\omega^n= \delta \omega^{n-1}$ which again induces the required identity since the differential in degree $-n$ for $\extAnti^{\bullet}$ is the inclusion of holomorphic $n$-forms, and this inclusion differential commutes with skew-symmetrization.
\item[\bf Step 2:\hspace{-2.5mm}] \hspace{1.3mm}
Next, we check compatibility with face maps. Indeed, the way that \eqref{EQ define BM cells} defines all $\ell$-cells by only using the charts $\rho_{i_0},\dots, \rho_{i_\ell}$ that come from the induced $\ell$-cell given by their faces in $\Chart(W)_k$ (and then symmetrizing these) ensures that the face maps are immediately respected, i.e., $d_j\circ \BMmap(W)_k=\BMmap(W)_{k-1}\circ d_j$. 

To check compatibility with the degeneracy $s_j$, note that $s_j$ on a chart $(\rho_0,\dots, \rho_\ell)$ repeats the $j^\text{th}$ chart $\rho_j$, and, thus, in the new charts, we have the composed maps $\phi_{0,1}, \dots, \phi_{j-1,j},id,\phi_{j,j+1},\dots, \phi_{\ell-1,\ell}$ with an additional identity map; cf. \eqref{EQU:phi_{j,j+1}-diagram}. However, the $\omega^q$, for $q\ge 1$, vanish when evaluated on such a tuple of maps containing the identity, and consequently   \eqref{EQ define BM cells} vanishes in this case as well. To see why $\omega^q$ vanishes, recall from the proof of \cite[theorem 10.9]{TT76} that $\omega^{\bullet} = g^{\bullet}\tau^{\bullet}$ is the cup product of two cochains, where $g^{\bullet}$ is defined by a recursive construction. Note that each $g^q$ is defined by applying a chain homotopy to an alternating sum of elements of the form of either $a^{r}g^{q-r}$ (if $r\ge 2$) or $\delta g^{q-1}$. By \cite[lemma 11.1]{TT76} the $a^q$ vanish when evaluated on a tuple containing the identity and it's straightforward to see why the \v{C}ech differential applied to $g^{q-1}$ also vanishes on such a tuple. The same is true for $\tau^{\bullet}$ by \cite[lemma 11.3]{TT76}, and so the claim follows. Thus, the only time we can get non-vanishing labels is when considering either the face map $d_j$ (in $\Chart(W)$) that combines $\phi_{j-1,j}\circ id$ or the face map $d_{j+1}$ that combines $id\circ \phi_{j,j+1}$, and for these two faces we get the same labels as the ones coming from the original chart $(\rho_0,\dots, \rho_\ell)$. According to equation \eqref{EQU:degen-sj(cI)}, this is equal to the application of the degeneracy $s_j$ in $\extAnti(W)=DK \left( {\extAnti^{\bullet}}(W)\right)$.

\item[\bf Step 3:\hspace{-2.5mm}] \hspace{1.3mm}
Finally, we show that the assignment is natural with respect to a biholomorphism $\psi: W' \to W$. Recall that for a $k$-simplex, $\{\rho_{i} : W \to V_i \}_{i=0}^k\in \Chart(W)_k$ as above, then $\Chart(\psi)(\rho_{\bullet}):= \{ \rho_{i}\circ \psi : W' \to (\rho_{i}\circ \psi) \left(V_i\right) \}_{i=0}^k\in \Chart(W')_k$. On the other hand for a $k$-simplex $x \in \extAnti(W)_k$, the $k$-simplex $\extAnti(\psi)(x) \in \extAnti(W')_k$ is given on each $\ell$-cell $(i_0, \ldots, i_{\ell})$ by $(\psi \times \psi)^*x_{i_0, \ldots, i_{\ell}}$. Then the claim is proven by: 
\begin{align*}
&\left(\extAnti(\psi)\left(\BMmap(W)_k(\rho_{\bullet}) \right)\right)_{i_0, \ldots, i_{\ell}}=(\psi \times \psi)^* \left( \BMmap(W)_k(\rho_{\bullet}) \right)_{i_0, \ldots, i_{\ell}} \\
&=\sum_{\sigma_\in S_{\ell+1}}\frac{ \text{sgn}(\sigma)}{(\ell +1)! }\cdot (\psi \times \psi)^*  \left(  \left(\rho_{i_{\sigma(\ell)}}\times \rho_{i_{\sigma(\ell)}} \right)^* \omega^{\ell}(\phi_{i_{\sigma(0)},i_{\sigma(1)}}, \ldots, \phi_{i_{\sigma(\ell -1)},i_{\sigma(\ell)}})  \right)\\
&=\sum_{\sigma_\in S_{\ell+1}}\frac{ \text{sgn}(\sigma)}{(\ell +1)! }\cdot \left( \left( \rho_{i_{\sigma(\ell)}}\circ \psi \right) \times \left( \rho_{i_{\sigma(\ell)}}\circ \psi \right) \right)^* \omega^{\ell}(\rho_{i_{\sigma(0)}} \circ \rho_{i_{\sigma(1)}}^{-1}, \ldots, \rho_{i_{\sigma(\ell-1)}} \circ \rho_{i_{\sigma(\ell)}}^{-1}) \\
&= \BMmap(W')_k \left(\Chart(\psi)\left( \rho_{\bullet} \right)  \right)_{i_0 , \ldots, i_{\ell}},
\end{align*}
where the last equality uses $\rho_{i_j} \circ \rho_{i_k}^{-1}   =\left(\rho_{i_j} \circ \psi\right) \circ \left( \rho_{i_k} \circ \psi\right)^{-1}$.
\end{itemize}
\end{proof}
\end{proposition}

\begin{definition}\label{DEF BM}
The \emph{Bochner-Martinelli $\delbar$-parametrix} is the morphism of $\infty$-prestacks $\BMmap: \Chart \to \extAnti$ defined in proposition \ref{PROP TTmap is a map}.
\end{definition}

Combining definitions \ref{DEF hartogs} and \ref{DEF BM}, we get the composition $\Hartmap\circ \BMmap: \Chart \to \extAnti\to \holFormP{n}$.

\section{Invariant polynomials}\label{SEC:Invariant-polynomials}

We now give another construction for a map of $\infty$-prestacks $\Cfmap{T}: \Chart \to \holFormP{k}$ associated to any fixed choice of a $GL$-invariant map $T$ on $k$ matrices. We review how such maps $T$ are given by sums and products of trace maps, and, when $T$ is symmetric, how this is related to invariant polynomials used in the theory of characteristic classes. The Bochner-Martinelli $\delbar$-parametrix, post-composed with the Hartogs map $ \Chart \xrightarrow{ \BMmap} \extAnti \xrightarrow{\Hartmap} \holFormP{n}$ turns out to be such a map of $\infty$-prestacks $\Hartmap \circ \BMmap=\Cfmap{Todd_n}$, where $T=Todd_n$ is the $n^{\text{th}}$ Todd symmetric polynomial, as was shown by Toledo and Tong in \cite{TT76}. This is our version of the Hirzebruch-Riemann-Roch theorem as stated in theorem \ref{THM main theorem}.

\begin{definition}
Let $k>0$. A $GL$-invariant map on $k$ many $n$-by-$n$ matrices is a linear map $T:\underbrace{\C^{n\times n}\ot\dots\ot \C^{n\times n}}_{k\text{ tensors}}\to \C$ such that
\[
T(A^{-1}\cdot M_1\cdot A,\dots,A^{-1}\cdot M_k\cdot A)=T(M_1,\dots, M_k), \quad \forall A\in GL(\C^n)
\]

For a choice of a permutation $\sigma\in S_k$, we can define such a map $T=T_\sigma$ as follows. If we write $\sigma$ as a list of cycles, $\sigma=(\sigma^{(1)}_1 \dots \sigma^{(1)}_{p_1})\dots (\sigma^{(r)}_1 \dots \sigma^{(r)}_{p_r})$, then we define the $GL$-invariant map $T_\sigma:(\C^{n\times n})^{\ot k}\to \C$ as a product of trace maps
\[
T_\sigma(M_1,\dots, M_k):=tr(M_{\sigma^{(1)}_1}\cdot  \ldots \cdot M_{\sigma^{(1)}_{p_1}})\cdot \ldots \cdot tr(M_{\sigma^{(r)}_1}\cdot  \ldots \cdot M_{\sigma^{(r)}_{p_r}})
\]
For example, if $\sigma=\left(\begin{matrix}1 & 2 & 3 & 4 & 5 & 6 & 7 & 8 \\ 3 & 1 & 7 & 4 & 8 & 5 & 2 & 6 \end{matrix}\right)=(1\,\, 3\,\, 7\,\, 2)(4)(5\,\, 8 \,\,6)$, then $T_\sigma(M_1,\dots, M_8)=tr(M_1 M_3 M_7 M_2)\cdot tr(M_4)\cdot tr(M_5 M_8 M_6)$, and clearly $T_\sigma$ is $GL$-invariant by the trace property $tr(M\cdot A)=tr(A\cdot M)$.

From \cite[theorem 5.10]{TT76} we see that any $GL$-invariant map $T$ is given by a linear combination of the $T_\sigma$, i.e., $T=\sum_{\sigma\in S_k} c_\sigma \cdot T_\sigma$ for some $c_\sigma\in \C$.
\end{definition}

For any $GL$-invariant map $T$, we obtain an induced map of $\infty$-prestacks from $\Chart$ to $\holFormP{k}$, i.e., a natural transformation of functors $\Cfmap{T}: \Chart \to \holFormP{k}$.

\begin{notation}\label{NOT:T[theta...]}
Recall from \eqref{EQU:phi_{j,j+1}-diagram}, that an $\ell$-simplex in $\Chart(W)$ is a sequence of maps $(\rho_0:W\to V_0,\dots, \rho_\ell:W\to V_\ell)$, with $V_j\subset \C^n$ for all $j$, and we denote, as usual, $\phi_{p,q}=\rho_p\circ \rho^{-1}_q:V_q\to V_p$. Following \cite[page 282]{TT76}, we set $\theta(\phi_{p,q}):=\dot{\phi}_{p,q}^{-1}\cdot d(\dot{\phi}_{p,q})$, where the ``dot'' denotes the Jacobian, so that $\theta(\phi_{p,q})\in \Omega^1(V_q)\ot \C^{n\times n}$ is a matrix-valued $1$-form.

The pullback under a biholomorphic map $\varphi:V_r\to V_q$ is given by $\varphi^\sharp:\Omega^1(V_q)\ot \C^{n\times n}\to \Omega^1(V_r)\ot \C^{n\times n}$, where $\varphi^\sharp(\alpha\ot M)(x)=(\varphi^*\alpha)(x)\ot (\dot{\varphi}^{-1}(x)\cdot M\cdot \dot{\varphi}(x))$. In particular, we define
\begin{equation}\label{EQU:theta_r}
\theta_r:=\phi_{r,\ell}^\sharp(\theta(\phi_{r-1,r})) \in \Omega^1(V_\ell)\ot \C^{n\times n}
\end{equation}

If $T$ is a $GL$-invariant map on $k$ matrices, we also denote by $T$ the induced map $T:(\Omega^1(V_q)\ot \C^{n\times n})^{\ot k}\to \Omega^k(V_q)$, $T[\alpha_1\ot M_1,\dots,\alpha_k\ot M_k]=\alpha_1\wedge\dots\wedge \alpha_k\cdot T(M_1,\dots, M_k)$. Note, that the $GL$-invariance of $T$ implies, that for the biholomorphic map $\varphi:V_r\to V_q$, we have
\begin{equation}\label{EQU:phi-sharp-pullback}
T[\varphi^\sharp(\alpha_1\ot M_1),\dots, \varphi^\sharp(\alpha_k\ot M_k)]
=\varphi^*(T[\alpha_1\ot M_1,\dots, \alpha_k\ot M_k])
\end{equation}
\end{notation}
Using the above notation, we will define a map of $\infty$-prestacks $\Cfmap{T}: \Chart \to \holFormP{k}$.
\begin{definition}\label{DEF cfmap is map}
Let $T$ be a $GL$-invariant map on $k$ matrices, and let $W$ be an object in $\Chart$. We define a map $\Cfmap{T}(W):\Chart(W)\to \holFormP{k}(W)$ of simplicial sets on $q$-simplicies $\Cfmap{T}(W)_q:\Chart(W)_q\to \holFormP{k}(W)_q$, i.e., for $\rho_{\bullet}= (\rho_0:W\to V_0,\dots, \rho_q:W\to V_q) \in \Chart(W)_q$ we assign a label $\Cfmap{T}(W)_q(\rho_\bu)_{i_0,\dots, i_\ell}$ for each $\ell$-cell given by some indices $0 \le i_0 < \cdots < i_{\ell} \le q$ as follows:
\begin{itemize}
\item
if $\ell=k$, we set $\wt\theta_r=\phi_{i_r,i_k}^\sharp(\theta(\phi_{i_{r-1},i_r}))$ for $r=1,\dots, k$, where as usual $\phi_{i,j}=\rho_i\circ \rho_j^{-1}$, and with this we define 
\begin{equation}
\Cfmap{T}(W)_q(\rho_0,\dots, \rho_q)_{i_0,\dots, i_k}:= \rho^*_{i_k}(T[\wt\theta_1,\dots, \wt\theta_k])\label{EQ Cfmap def}
\end{equation}
\item
if $\ell\neq k$, we set $(\Cfmap{T}(W)_q(\rho_\bu))_{i_0,\dots, i_\ell}:=0$
\end{itemize}
\end{definition}
\begin{proposition}\label{PROP cfmap is map}
The above definition \ref{DEF cfmap is map} gives a map of $\infty$-prestacks $\Cfmap{T}: \Chart \to \holFormP{k}$.
\end{proposition}
\begin{proof}
We need to show that $\Cfmap{T}$ is a natural transformation. We do this in three steps. First, we show that the above assignment gives a well-defined $\ell$-simplex, second, the assignment $\Cfmap{T}(W)$ respects face and degeneracy maps, and third, $\Cfmap{T}$ gives the correct commutative diagram for a natural transformation.
\begin{itemize}
\item[\bf Step 1:\hspace{-2.5mm}] \hspace{1.3mm}
To see that the above are well-defined $\ell$-cells in $\holFormP{k}(W)_\ell$, we need to check that the alternating sum of labelings of the faces of any $(k+1)$-cell vanishes. It is enough to show this for a fixed $(k+1)$-simplex $(\rho_0,\dots, \rho_{k+1})$, for which the sum of the faces becomes:
\begin{multline}\label{EQU:sum-faces-theta}
\rho^*_{k+1}\left(T\left[\theta_2,\dots, \theta_{k+1}\right]\right)+\sum_{j=1}^k(-1)^j \rho^*_{k+1}\left(T\left[\theta_1,\dots,\phi^\sharp_{j+1,k+1}\left(\theta\left(\phi_{j-1,j}\circ \phi_{j,j+1}\right)\right) ,\dots,\theta_{k+1}\right]\right)
\\+(-1)^{k+1} \rho^*_{k}
\left(T\left[\phi^\sharp_{1,k}\left(\theta\left(\phi_{0,1}\right)\right),\dots, \phi^\sharp_{k,k}\left(\theta\left(\phi_{k-1,k}\right)\right)\right]\right)
\end{multline}
Now, $\theta$ applied to a composition is given by 
\begin{equation}
\theta(g\circ h)=h^\sharp(\theta(g))+\theta(h),\label{EQ theta g circ h}
\end{equation}
since
\begin{multline*}
(({g\circ h})\dot{\,\,})^{-1}\cdot d({g\circ h})\dot{} = (\dot{g}|_{h}\cdot \dot{h})^{-1}\cdot d(\dot{g}|_{h}\cdot \dot{h})
= \dot{h}^{-1}\cdot \dot{g}|_{h}^{-1}\cdot \Big(d(\dot{g}|_{h})\cdot \dot{h}+  \dot{g}|_{h}\cdot d(\dot{h})\Big)
\\
= \dot{h}^{-1}\cdot \dot{g}|_{h}^{-1}\cdot d(\dot{g}|_{h})\cdot \dot{h}+  \dot{h}^{-1}\cdot d(\dot{h})
= \dot{h}^{-1}\cdot h^*(\dot{g}^{-1} d\dot{g})\cdot \dot{h}+  \theta(h)=h^\sharp(\theta(g))+\theta(h)
\end{multline*}
In the sum over $j$ in \eqref{EQU:sum-faces-theta}, we simplify the composition of $\phi$s as
\begin{equation}\label{EQU:phi^sharp(thetaotheta)}
\phi^\sharp_{j+1,k+1}(\theta(\phi_{j-1,j}\circ \phi_{j,j+1}))
\substack{\eqref{EQ theta g circ h}\\=}\underbrace{\phi^\sharp_{j+1,k+1}\circ \phi_{j,j+1}^\sharp}_{=(\phi_{j,j+1}\circ \phi_{j+1,k+1})^\sharp}(\theta(\phi_{j-1,j}))+\phi^\sharp_{j+1,k+1}(\theta( \phi_{j,j+1}))
=\theta_j+\theta_{j+1}
\end{equation}
so that \eqref{EQU:sum-faces-theta} becomes a telescopic sum, since, the expression on the left-hand side of \eqref{EQU:phi^sharp(thetaotheta)} appears in the $j^\text{th}$ term of the sum in \eqref{EQU:sum-faces-theta}, and we may replace it with $\theta_j+\theta_{j+1}$. This $j^\text{th}$ summand then becomes $(-1)^j \cdot\big(\rho^*_{k+1}\big(T\big[\dots, \widehat{\theta_{j+1}}, \dots\big]\big)+\rho^*_{k+1}\big(T\big[\dots, \widehat{\theta_j}, \dots\big]\big)\big)$, and so when we take the sum over $j$ we ``skip'' over any $\theta_i$ for $i=2,\dots,k$ twice with opposite signs, leaving only a first and last term. Now, the first term, which skips over $\theta_1$, cancels with the very first term in \eqref{EQU:sum-faces-theta} (before the big sum over $j$), and similarly the last term, which skips over $\theta_{k+1}$, cancels with the very last term in \eqref{EQU:sum-faces-theta}, since $\rho_k= \phi_{k,k+1}\circ \rho_{k+1}$ and using \eqref{EQU:phi-sharp-pullback}, we have that
$$
 \rho^*_{k}(T[\phi^\sharp_{1,k}(\theta(\phi_{0,1})),\dots, \phi^\sharp_{k,k}(\theta(\phi_{k-1,k}))])
=
 \rho^*_{k+1}(T[\phi^\sharp_{1,k+1}(\theta(\phi_{0,1})),\dots, \phi^\sharp_{k,k+1}(\theta(\phi_{k-1,k}))])
$$
Thus, the expression in \eqref{EQU:sum-faces-theta} vanishes, which is what we needed to show.
\item[\bf Step 2:\hspace{-2.5mm}] \hspace{1.3mm}
To check that $\Cfmap{T}(W):\Chart(W)\to \holFormP{k}(W)$ is a map of simplicial sets, we need to check that it respects the face and degeneracy maps. Compatibility with face maps follows by construction, since $k$-cells in $\holFormP{k}(W)_q$ are labeled by applying $\Cfmap{T}(W)$ to the corresponding faces in $\Chart(W)_q$.

Now, for a degeneracy $s_j$ applied to a chart $(\rho_0,\dots, \rho_k)$, we get the transition maps $\phi_{0,1}, \dots, \phi_{j-1,j},id,\phi_{j,j+1},\dots, \phi_{k-1,k}$ with one added identity map, as described in step 2 of the proof of proposition \ref{PROP TTmap is a map}. However, $\theta(id)=0$, and so, the only non-vanishing labels are given when we apply the face map $d_j$ or $d_{j+1}$ in $\Chart(W)$, and on these two faces we get labels coming from the original chart $(\rho_0,\dots, \rho_k)$. By equation \eqref{EQU:degen-sj(cI)}, this is equal to the application of the degeneracy $s_j$ in $\holFormP{k}(W)=DK(\holForm^k(W)[-k])$.

\item[\bf Step 3:\hspace{-2.5mm}] \hspace{1.3mm}
If $\phi:W'\to W$ is a morphism in $\Bihol$, i.e., a biholomorphism onto its image, then the diagram for a natural transformation commutes
\[
\begin{tikzcd}
\Chart(W) \arrow[swap]{d}{\Cfmap{T}(W)}\arrow{rr}{\Chart(\phi)}&& \Chart(W') \arrow{d}{\Cfmap{T}(W')}\\
\holFormP{k}(W) \arrow{rr}{\holFormP{k}(\phi)} && \holFormP{k}(W')
\end{tikzcd}
\]
since both $\Chart(\phi)$ and $\holFormP{k}(\phi)$ are given by pullback by $\phi$, and the induced charts under $\phi$ are just $\rho_k|_{\phi(W)}\circ \phi$ giving the induced ``$\phi_{p,q}$'' maps as restrictions of the previous maps $(\rho_p|_{\phi(W)}\circ \phi)\circ (\rho_q|_{\phi(W)}\circ \phi)^{-1}=\phi_{p,q}|_{Im(\rho_q|_{\phi(W)}\circ \phi)}$. This shows that above diagram commutes, and thus $\Cfmap{T}$ is a map of $\infty$-prestacks.
\end{itemize}
\end{proof}

Many interesting cases of $GL$-invariant maps have an additional property of being symmetric in their inputs of the $k$ matrices.

\begin{example}\quad
\begin{enumerate}
\item
If $\sigma, \tau\in S_k$ are two permutations on $k$ elements, and $\sigma$ is given by the list of cycles 
$
\sigma=(\sigma^{(1)}_1 \dots \sigma^{(1)}_{p_1})\dots (\sigma^{(r)}_1 \dots \sigma^{(r)}_{p_r}),
$
then $\tau\circ \sigma\circ \tau^{-1}$ is given by the list of cycles
\[
\tau\circ \sigma\circ \tau^{-1}=\Big(\tau(\sigma^{(1)}_1) \dots \tau(\sigma^{(1)}_{p_1})\Big)\dots \Big(\tau(\sigma^{(r)}_1) \dots \tau(\sigma^{(r)}_{p_r})\Big)
\]
In particular, for the transposition $\tau=(ij)$ that switches $i$ and $j$, we get that\linebreak $T_{\sigma}(M_1,\dots,M_i,\dots, M_j, \dots, M_k)=T_{\tau\circ \sigma\circ \tau^{-1}}(M_1,\dots,M_j,\dots, M_i, \dots, M_k)$. Thus, we can define a symmetrized version of $T_\sigma$ by setting
\[
T^{Sym}_{\sigma}:=\frac 1 {k!}\sum_{\tau\in S_k} T_{\tau\circ \sigma\circ \tau^{-1}}
\]

Since each $T^{Sym}_{\sigma}$ is symmetric with respect to its inputs $M_1,\dots, M_k$, it is given by symmetrized trace maps
\begin{equation}\label{EQU:T^Sym_sigma}
\frac 1 {k!} \sum_{\tau\in S_k} tr(M_{\tau(1)}\dots M_{\tau(p_1)})\cdot tr(M_{\tau(p_1+1)}\dots M_{\tau(p_2)})\cdot  \ldots \cdot tr(M_{\tau(p_1+\dots +p_{r-1}+1)}\dots M_{\tau(p_r)}),
\end{equation}
thus it only depends on the numbers of inputs $p_1,\dots, p_r$ in each trace. 
\item
From Chern-Weil theory, we know that the invariant polynomials on $gl(n,\C)$ are generated as an algebra by the trace functions $\T_j(M)=tr(M^j)$ for $M\in gl(n,\C)$, or alternatively by a \emph{symmetric polynomial} in the eigenvalues of $M$, for example, $\T_j(M)=tr(M^j)=t_1^j+\dots + t_n^j$, where $t_i$ are the eigenvalues of $M$ with possible repetition due to their multiplicies; see e.g. \cite[Appendix B]{Tu}.

Note that to a product $\T_{p_1}\cdot \ldots\cdot \T_{p_r}$ we can associate the $GL$-invariant map $T^{Sym}_\sigma$ from \eqref{EQU:T^Sym_sigma} for a permutation $\sigma$ with $r$ cycles of length $p_1,\dots, p_r$, respectively.
\item
For any symmetric polynomial on the eigenvalues $t_1,\dots, t_n$ of a matrix $M$, we get an induced $GL$-invariant polynomial. For example, the polynomial associated to the Chern character is the symmetric polynomial 
\begin{align*}
\mc Ch(M)=&e^{t_1}+\dots+e^{t_n}=n+(t_1+\dots +t_n)+\frac 1 {2!} (t_1^2+\dots +t_n^2)+\dots
\\
=&n+\mc T_1(M)+\frac 1 {2!}\mc T_2(M)+\frac 1 {3!} \mc T_3(M)+\dots
\end{align*}
In particular, the $k^{\text{th}}$ Chern character polynomial is $\mc Ch_k=\frac 1 {k!}\T_k$.
\item
An alternative set of generators for invariant polynomials on $gl(n,\C)$ is given by the elementary symmetric polynomials $\mc S_j(M)=\sum\limits_{1\leq i_1<i_2<\dots<i_j\leq n} t_{i_1}\cdot t_{i_2}\cdot \ldots \cdot t_{i_j}$ where again $t_1,\dots, t_n$ are the eigenvalues of $M$. We also call $\mc S_k$ the $k^{\text{th}}$ Chern class polynomial.

The elementary symmetric polynomials and the trace polynomials are related via Newton's identities:
\[
\mc T_j-\mc T_{j-1}\mc S_1+\mc T_{j-2}\mc S_2-\dots +(-1)^{j-1}\cdot \mc T_{1}\mc S_{j-1}+(-1)^j\cdot j\cdot \mc S_j=0
\]
\item
Another example is given by the symmetric Todd polynomial in the eigenvalues $t_1,\dots, t_n$ of $M$:
$$
\mc Todd(M)=\frac{t_1}{1-e^{-t_1}}\cdot \frac{t_2}{1-e^{-t_2}}\cdot \ldots \cdot \frac{t_n}{1-e^{-t_n}}
$$
The $k^\text{th}$ (homogeneous) component of the Todd polynomial is denoted by $\mc Todd_k$. Expanding the above via the Taylor series $\frac t {1-e^{-t}}=1+\frac t 2+\frac{t^2}{12}-\frac{t^4}{120}+\dots$, we get in low dimensions that
\begin{align*}
\mc Todd_1(M)=&\frac 1 2 \mc S_1(M) = \frac 1 2 \mc T_1(M)
\\
\mc Todd_2(M)=&\frac 1 {12}( \mc S_1(M)^2+\mc S_2(M))=\frac 1 {24}(3\cdot \mc T_1(M)^2-\mc T_2(M))
\\
\mc Todd_3(M)=&\frac 1 {24} \mc S_1(M)\cdot \mc S_2(M)=\frac 1 {48}( \mc T_1(M)^3-\mc T_1(M)\cdot \mc T_2(M))
\end{align*}
We denote by $Todd_k$ the $GL$-invariant map corresponding to $\mc Todd_k$.
\end{enumerate}
\end{example}

The following theorem is a reinterpretation of \cite[corollary 11.5]{TT76}.
\begin{theorem}\label{THM main theorem}
The composition of the Hartogs map (definition \ref{DEF hartogs}) and the Bochner-Martinelli $\delbar$-parametrix (definition \ref{DEF BM}) $ \Chart \xrightarrow{\BMmap} \extAnti \xrightarrow{\Hartmap} \holFormP{n}$ is equal to $\Cfmap{Todd_n}$ as maps of $\infty$-stacks:
\begin{equation}\label{EQU:HartoBM=CToddn}
\Hartmap \circ \BMmap=\Cfmap{Todd_n},
\end{equation}
\end{theorem}
\begin{proof}
The construction of the three maps of $\infty$-prestacks involved is the content of propositions \ref{PROP Hartmap is a map}, \ref{PROP TTmap is a map}, and \ref{PROP cfmap is map}. All that remains to check is that the two natural transformations in \eqref{EQU:HartoBM=CToddn} agree on objects, which we now show to follow from the Universal Riemann-Roch Theorem (URRT)  ``$\wt\kappa=\mathcal T_n(\theta^n)$'' stated in \cite[corollary 11.5]{TT76}.

Consider an object $W$ of $\Bihol$. Then, by definition \ref{DEF cfmap is map}, $\Cfmap{Todd_n}(W)$ is non-zero only on $n$-cells, where it is given by $\Cfmap{Todd_n}(W)_n(\rho_0,\dots, \rho_n)_{0,\dots, n}=\rho^*_n(Todd_n[\theta_1,\dots, \theta_n])$, cf. \eqref{EQ Cfmap def}, and similarly for $n$-cells that are faces of higher dimensional simplices. Here, our notation \ref{NOT:T[theta...]} follows the one from \cite[second paragraph on page 282]{TT76} showing that in the notation of \cite{TT76}, this equals $f(\theta^n)$ applied to the transition functions $\phi_{0,1},\dots,\phi_{n-1,n}$, which is then pulled back to $W$, and where $f=Todd_n$ is the $n^\text{th}$ Todd polynomial. This is the right-hand side ``$\mathcal T_n(\theta^n)$''  of the URRT in \cite[corollary 11.5]{TT76} (before the pullback to $W$).

Looking at the left-hand side of \eqref{EQU:HartoBM=CToddn}, by \eqref{EQ Hartog Delta} and example \ref{DEF:DK(V)}, $\Hartmap \circ \BMmap(W)$ is non-zero only on $n$-cells, where it is given by $\Hartmap \circ \BMmap(W)_n(\rho_0,\dots,\rho_n)_{0,\dots,n}=\Delta^*(\overline{\omega} )$, with $\omega=\sum_{\sigma\in S_{n+1}} \frac{\text{sgn}(\sigma)}{(n +1)! }\cdot \left( \rho_{\sigma(n)} \times \rho_{\sigma(n)}\right)^*\left(\omega^{n}(\phi_{\sigma(0),\sigma(1)}, \ldots, \phi_{\sigma(n-1),\sigma(n)}) \right)$ according to \eqref{EQ define BM cells}. Thus the holomorphic extension and pullback to the diagonal $\Delta$ of $\omega$ is the skew-symmetrized version of the ``$\kappa$'' term from \cite[page 282]{TT76} coming from the skew symmetrized \v{C}ech parametrix, which is then pulled back to $W$. According to \cite[theorem III]{TT76}, the skew-symmetrized term ``$\wt\kappa$'' (before pullback to $W$),  is the term on the left-hand side of the URRT stated in \cite[corollary 11.5]{TT76}.

The claim of the theorem thus follows from the URRT \cite[corollary 11.5]{TT76}  pulled back to $W$. Since the theorem is a just a reinterpretation of \cite[corollary 11.5]{TT76}, it may be worth to briefly outline the steps for its proof in \cite{TT76}.
\begin{enumerate}[(i)]
\item\label{item:11.5-i}
Starting from the Bochner-Martinelli kernel $\omega^0$ in \cite[(2.10)]{TT76}, the higher terms in the parametrix $\omega^1,\dots, \omega^{n-1}$ with $\delbar \omega^q=\delta\omega^{q-1}$ are constructed in \cite[chapters 7-10]{TT76}, with their final result stated in \cite[theorem 10.9]{TT76}. Moreover $\omega^n$ is defined as $\delta \omega^{n-1}$, see \cite[page 282]{TT76}.
\item\label{item:11.5-ii}
A sufficient condition for when a cocyle (such as the holomorphic extension $\Delta^*\omega^n$ where $\omega^n$ is from \ref{item:11.5-i}) is of the form $f(\theta^n)$ for some $GL$-invariant map $f$ is stated in \cite[theorem II]{TT76}. \cite[theorem II]{TT76} is proved in \cite[chapter 5]{TT76}.
\item\label{item:11.5-iii}
\cite[theorem III]{TT76} states that the conditions for \cite[theorem II]{TT76} hold for the skew-symmetrized holomorphic extension $\wt\kappa:=\Delta^*\wt\omega^n$ where, again, $\omega^n$ is from \ref{item:11.5-i}. \cite[theorem III]{TT76} is proved in \cite[chapter 11]{TT76}.
\item\label{item:11.5-iv}
Finally, the URRT $\wt\kappa=\Delta^*\wt\omega^n=Todd_n(\theta^n)$ \cite[corollary 11.5]{TT76} follows by identifying the $GL$-invariant polynomial $f$ from \ref{item:11.5-ii} via multiplicative properties of the Chern numbers corresponding to $f$, see \cite[third paragraph on page 282]{TT76}.
\end{enumerate}
We also note that \cite[theorems II\ensuremath{'} and III\ensuremath{'} and corollary 12.1]{TT76} contain generalized versions of these statements concerning the analogue problem with coefficients in a bundle.
\end{proof}

Since the local projective model structure $\sPre(\Bihol)_{/{\mc S}}$ (definition \ref{DEF loc proj mod}) is a simplicial model category then its \emph{derived mapping space}\footnote{See for example \cite[section 17]{Hi} where there is a complete treatment via homotopy function complexes.} is given by
\begin{equation}\label{EQ RHom}
\mathbb{R}\text{Hom}(\mc F,\mc G) := \underline{\sPre}(\Bihol) \left(\check{\mc F} ,{\mc G}^{\dagger}\right).
\end{equation}
where $\check{\mc F}$ and ${\mc G}^{\dagger}$ are respectively a cofibrant replacement of $\mc F$ and a fibrant replacement of $\mc G$ with respect to the local projective module structure. Naturality of the (derived) mapping spaces means that the equality of maps of $\infty$-prestacks from theorem \ref{THM main theorem} immediately provides equality of maps of hom-spaces.
\begin{corollary}\label{COR main corollary}
For any $\infty$-prestack $\mc F$ on $\Bihol$ the pair of parallel morphisms:
\begin{equation}\label{EQ BM=Todd spre}
\begin{tikzcd}[column sep=huge]
\mathbb{R}\text{Hom}( \mc F, \Chart)
  \arrow[bend left=15]{r}[name=U, label=above:$\left(\Hartmap \circ \BMmap\right)_*$]{}
  \arrow[bend right=15]{r}[name=D,label=below:${\Cfmap{Todd_n}}_*$]{} &
\mathbb{R}\text{Hom}( \mc F, \holFormP{n})
  \arrow[shorten <=20pt,shorten >=15pt,equal,to path={(U) -- (D)}]{}.
\end{tikzcd}
\end{equation}
are equal as maps of simplicial sets.
\end{corollary}

\section{Evaluating $\infty$-stacks on $G$-Manifolds}\label{SEC eval on man}
Given an $n$-dimensional complex manifold, $M$, define the $\infty$-prestack $\varepsilon M$ on $\Bihol$ by assigning to each object the discrete simplicial set $\varepsilon M(W)$ of holomorphic embeddings $W \to M$. This $\infty$-prestack is not necessarily cofibrant in $\sPre(\Bihol)$ and neither is the analogous presheaf $\mmod{M}{G}$ (definition \ref{DEF M mmod G}) associated to a $G$-action on $M$, and so in those cases a cofibrant replacement is required in order to compute the derived mapping space \eqref{EQ RHom} as explained below. 
\begin{definition}\label{DEF G man}
Let $G$ be a discrete group. A \emph{$G$-manifold} is an $n$-dimensional complex manifold, $M$ with (bi-)holomorphic group action by discrete group $G$, ie $\rho: G \to Aut(M)$ is a group homomorphism where, for each $g \in G$, $\rho(g)$ provides a biholomorphism of $M$. 
\end{definition}
The following natural definition adapts the nerve of the action groupoid to our setting of $\infty$-prestacks:
\begin{definition}\label{DEF M mmod G}
Let $M$ be a $G$-manifold. Define the \emph{quotient $\infty$-prestack} $\mmod{M}{G}$ to be the following $\infty$-prestack on $\Bihol$:
\begin{itemize}
\item The presheaf of $n$-simplices is  $\left(\mmod{M}{G} \right)_n:= \coprod\limits_{g_1, \ldots, g_n} \varepsilon M$,
\item The simplicial face maps $d_k: \left(\mmod{M}{G} \right)_{n+1} \to \left( \mmod{M}{G} \right)_n$ are given object-wise by 
\[\left( \left(g_1, \ldots, g_{n+1}\right), W \xrightarrow{f} M\right) \xmapsto{d_k} \begin{cases}
\left(\left(g_2, g_{3}, \ldots  g_{n+1}\right), W \xrightarrow{f} M\right) & \text{for} \quad k=0\\
\left(\left(g_1, \ldots, g_{k} \cdot g_{k+1}, \ldots  g_{n+1}\right), W \xrightarrow{f} M\right)& \text{for} \quad 0<k<n+1\\
\left(\left(g_1,  \ldots , g_{n}\right), W \xrightarrow{\rho(g_{n+1}) \circ f} M\right)& \text{for} \quad k=n+1
\end{cases}  \] 
\item The simplicial degeneracy maps $s_k: \left(\mmod{M}{G} \right)_{n} \to \left(\mmod{M}{G} \right)_{n+1}$ are given object-wise by 
\[\left( \left(g_1, \ldots, g_{n}\right), W \xrightarrow{f} M\right) \xmapsto{s_k} \left(\left( \ldots,g_{k-1}, 1,g_k, \ldots \right), W \xrightarrow{f} M\right) .  \] 
\end{itemize}
\end{definition}

\begin{definition}\label{DEF equiv good cover}
An open cover $\mathcal{U}= \{U_i \subset M\}_{i \in I}$ of a $G$-manifold $M$, with action $\rho: G \times M \to M$, is called a \emph{$G$-equivariant $\Bihol$-cover} if each $U_i$ is an object in $\Bihol$ and $G$ acts on the indices $I \times G \to I$ written $(i,g) \mapsto i \cdot g$ so that, for each $g$, $\rho(g): U_i \to U_{i \cdot g^{-1}}$ is a biholomorphism between open sets in $\mathcal{U}$. 
\end{definition}
Note that in the case when $G = \{ *\}$ and thus $\mmod{M}{G} = \varepsilon M$, a $G$-equivariant $\Bihol$-cover is exactly an open cover of $M$ by objects in $\Bihol$. In what follows, notation analogous to definition \ref{DEF covering family} is used but here the situation is simpler since the $U_i$'s are open subset of a fixed manifold so $U_{i_0, \ldots, i_p}$ is just defined as the intersection $U_{i_0} \cap \ldots \cap U_{i_p}$. 
\begin{proposition}\label{PROP bi cover pre}
Given a $G$-equivariant $\Bihol$-cover $\mathcal{U}$ of $M$, the following data forms a bi-simplicial presheaf $\mmod{\mathcal{U}}{G}$ on $\Bihol$:
\begin{itemize}
\item The bi-simplicial presheaf is defined bi-level-wise by  \[\left(\mmod{\mathcal{U}}{G} \right)_{m,n}:= \coprod\limits_{\substack{ i_0, \ldots, i_m\\ g_1, \ldots, g_n} } {\yo}U_{\left(i_0,\ldots, i_m\right) \cdot \left(g_1 \cdots g_n\right)},\] where for each $g \in G$, the above multiplication on tuples in $I$ is defined $U_{\left(i_0,\ldots, i_m\right) \cdot g}: = U_{i_0\cdot g,\ldots, i_m \cdot g }$,
\item The face maps {$d^{-1,0}_{k} \left(\mmod{\mathcal{U}}{G} \right)_{m+1,n} \to\left(\mmod{\mathcal{U}}{G} \right)_{m,n}$} are given object-wise by 
\begin{multline*}
{\hspace{1cm}\left( \left( \substack{i_0, \ldots, i_{m+1}\\ g_1, \ldots , g_n}\right), W \xrightarrow{f} U_{\left(i_0, \ldots, i_{m+1}\right) \cdot  \left(g_1, \ldots, g_{n}\right)}\right)}
\\ {\xmapsto{d^{-1,0}_k} 
\left( \left( \substack{i_0, \ldots, \widehat{i_k} , \ldots , i_{m+1}\\ g_1, \ldots , g_n}\right),W \xrightarrow{f} U_{\left(i_0, \ldots , i_{m+1}\right) \cdot  \left(g_1 \dots g_{n}\right)}\hookrightarrow U_{\left(i_0, \ldots, \widehat{i_k} , \ldots , i_{m+1}\right) \cdot \left(g_1 \cdots g_{n}\right)  }    \right)}
\end{multline*}
\item The face maps {$d^{0,-1}_{k}: \left(\mmod{\mathcal{U}}{G} \right)_{m,n+1} \to\left(\mmod{\mathcal{U}}{G} \right)_{m,n}$} are given object-wise by 
\begin{align*}
&\left( \left( \substack{i_0, \ldots, i_m\\ g_1, \ldots , g_{n+1}}\right), W \xrightarrow{f} U_{\left(i_0, \ldots, i_m\right) \cdot  \left(g_1 \cdots g_{n+1}\right)} \right) \\
{\xmapsto{d^{0,-1}_k}} &\begin{cases}
\left( \left( \substack{i_0\cdot g_1, \ldots, i_m \cdot g_1\\ g_2, \ldots , g_{n+1}}\right), W \xrightarrow{f} U_{\left(i_0\cdot g_1, \ldots, i_m\cdot g_1\right) \cdot  \left(g_2 \cdots g_{n+1}\right)}\right) & \text{for} \quad k=0\\
\left( \left( \substack{i_0, \ldots, i_m\\ g_1, \ldots, g_k \cdot g_{k+1} , \ldots, g_{n+1}}\right), W \xrightarrow{f} U_{\left(i_0, \ldots, i_m\right) \cdot  \left(g_1\cdots \left(g_k \cdot g_{k+1}\right) \cdots g_{n+1}\right)}\right)& \text{for} \quad {0<k<n+1}\\
\left( \left( \substack{i_0, \ldots, i_m\\ g_1, \ldots , g_{n}}\right), W \xrightarrow{f} U_{\left(i_0, \ldots, i_m\right) \cdot  \left(g_1 \cdots g_{n+1} \right)}\xrightarrow{{\rho(g_{n+1})}} U_{\left(i_0, \ldots, i_m\right) \cdot  \left(g_1 \cdots g_{n} \right)} \right)& \text{for} \quad {k=n+1}
\end{cases} 
\end{align*}
{where $\rho(g_{n+1})$ maps $U_{i\cdot g_{n+1}}\to U_i$ for any index $i$.}
\item The degneracy maps are appropriately given by either inserting a $1 \in G$ into a tuple of elements in $G$ or repeating an index $i_k$ in a tuple of indices from $I$.
\end{itemize}
We define the $\infty$-prestack $\abs{ \mmod{\mathcal{U}}{G}}$ to be the object-wise diagonal, or geometric realization, of this bi-simplicial presheaf. 
\end{proposition}

Alternatively, $\abs{\mmod{\mathcal{U}}{G}}$ could be defined as the nerve of a presheaf of groupoids. Given a test object $W$, consider the following category:
\begin{itemize}
\item An object $(i,f)$ is a morphism $W \xrightarrow{f} U_i$ for some $i \in I$. 
\item Given a pair $i,j \in I$, if $W \xrightarrow{f_i} U_i$ does not factor through $U_{j\cdot g,i}$ for some $g$ then the set of morphisms from $(j,f_j)$ to $(i, f_i)$ is empty. Otherwise, a morphism $(j, f_j) \to (i, f_{i})$ is a group element $g \in G$ so that $f_i$ factors through $U_{j \cdot g, i}$ and the composition 
\begin{equation} W \xrightarrow{f_i} U_{j \cdot g, i} \xrightarrow{{\rho(g)}} U_{j , i\cdot g^{-1}} \hookrightarrow U_j \label{EQ edge U mmod G gpd}\end{equation} is equal to $f_j$, where the unlabeled inclusion induces the face map $d_1$ from edges to vertices in $\abs{\mmod{\mathcal{U}}{G}}(W)$. 
\end{itemize}
Once the simplices of the nerve of the above category are unravelled, the equivalence with $\abs{\mmod{\mathcal{U}}{G}}(W)$ becomes clear. For example, consider the 2-simplex, $(k, f_k) \xrightarrow{h} (j, f_j) \xrightarrow{g} (i, f_{i})$. Applying \eqref{EQ edge U mmod G gpd} to both the domain-codomain matching conditions forced at each vertex, as well as the commuative triangle, the core data of such a $2$-simplex forms the tetrahedron:
\begin{equation}
\begin{tikzcd}
& W \arrow[dl, "f_i"'] \arrow[dr, "f_k"] & \\
U_{k \cdot (gh), j \cdot g, i} \arrow[rr, "\rho(hg)"', near start] \arrow[dr, "\rho(g)"'] & & U_{k, j \cdot h^{-1}, i \cdot (hg)^{-1}} \\
& U_{k \cdot h, j, i \cdot g^{-1}} \arrow[ur, "\rho(h)"'] 
\arrow[from=1-2, crossing over, "f_j" near start] & 
\end{tikzcd}
\end{equation}
One can check that this data is equivalent to the $2$-simplex $\left(g_1 g_2 , W \xrightarrow{f} U_{(i_0, i_1, i_2) \cdot (g_1 g_2)} \right)$ in $\abs{\mmod{\mathcal{U}}{G}}(W)_2$ after setting $g_1 := h$, $g_2 := g$, $i_0:= k$, $i_1:= j \cdot h^{-1}$, and $i_2 :=  i \cdot (hg)^{-1}$.

The consequences of this equivalent construction are recorded in the following lemma.
\begin{lemma}\label{LEM UmmodG 2cosk}
The $\infty$-prestack $\abs{\mmod{\mathcal{U}}{G}}$ is the nerve of a presheaf of groupoids and so takes values in $2$-coskeletal Kan complexes. \end{lemma}
In particular, since ${\mmod{M}{G}}$ is equal to $\abs{\mmod{\mathcal{U}}{G}}$ when $\mc U=\{M\subset M\}$ is the singleton-identity cover, and  $\Cech\mathcal{U}$ is equal to $\abs{\mmod{\mathcal{U}}{G}}$ when $G$ is a point, the same is statement applies to both of these $\infty$-prestacks.

\begin{lemma}\label{LEM iso hmtpy grps}
The natural map, $r: \abs{ \mmod{\mathcal{U}}{G}}\to \mmod{M}{G}$, induces maps of presheaves of sets $\pi_i(r): \pi_i\left(\abs{ \mmod{\mathcal{U}}{G}} \right)\to \pi_i\left(\mmod{M}{G}\right)$ which are isomorphisms after sheafification for $i=0,1$, provided each $U_i$ in $\mathcal{U}$ is connected.
\begin{proof}
Fix a test object $W \in \Bihol$. The map $r$ assigns on $p$-simplices:
\[ \left( W \xrightarrow{a} U_{(i_0,\ldots,  i_p)\cdot (g_1, \ldots,  g_p)} \right)\xmapsto{r} \left( W \xrightarrow{a} U_{(i_0,\ldots,  i_p)\cdot (g_1 , \ldots,  g_p)}  \hookrightarrow M \times \{ (g_1 , \ldots,  g_p)\}\right),\]
where $M \times \{ (g_1 , \ldots,  g_p)\}$ is representing the corresponding component of $\varepsilon M(W)$. This map is injective on $\pi_0$ before sheafification: if the images of two vertices $W \xrightarrow{a} U_i$ and $W \xrightarrow{a'} U_{i'}$ are related in $\varepsilon M(W)$, then there exists a map $b: W \to M \times{g}$ where $r(a) = b$ and $r(a') = \rho(g) \circ b$, but this map induces an edge $W \xrightarrow{a'} U_{(i, i' \cdot g)\cdot (g^{-1})}$ relating $a$ and $a'$ back in  $\abs{ \mmod{\mathcal{U}}{G}}(W)$. Surjectivity on $\pi_0$ requires sheafification (see \cite[lemma 6.8]{DI})  but amounts to: if $W \xrightarrow{b} M$ is a point in $ \mmod{M}{G}(W)$ then by covering $W$ with $\mathcal{W} = \{ W_i := b^{-1}(U_i)\}$ and taking the corresponding sieve $\overline{\mathcal{W}}$, for any $X \in \overline{\mathcal{W}}$, the restriction of $b$ to $X$ factors through some $U_i$ and so is locally in the image of the corresponding vertex from $\abs{ \mmod{\mathcal{U}}{G}}(W)$.

Next, we show injectivity on $\pi_1$ before sheafification (which preserves injectivity): the kernel of the map is trivial since the only null-homotopic loops in $\mmod{M}{G}(W)$ at a basepoint $W \xrightarrow{f} M$ are given by $e \in G$ and so the only loops in $\abs{ \mmod{\mathcal{U}}{G}}(W)$ which get sent to those loops are null-homotopic themselves. Finally, we show surjectivity on $\pi_1$ after sheafification: given a loop in $\mmod{M}{G}(W)$ at a basepoint $W \xrightarrow{f} M$ presented by $g \in G$, then $g$ fixes the map $f$. By holomorphicity, $g$ fixes the connected component of $M$ containing the image of $f$, and so it acts by the identity on every $U_i$ intersecting the image of $f$ (by connectedness). Thus covering $W$ by $W_i = f^{-1}(U_i \cap f(W))$ we obtain restriction basepoints $W_i \xrightarrow{f_i} U_i$ in $\abs{ \mmod{\mathcal{U}}{G}}(W_i)$ and loops presented by $g$, which each map to the corresponding loop in $\mmod{M}{G}(W_i)$.
\end{proof}
\end{lemma}

\begin{proposition}\label{PROP cofib repl M mmod G}
Let $M$ be a $n$-dimensional $G$-manifold and $\mathcal{U}:= \{ U_i \subset M \vert U_i \in \Bihol \}$ be an equivariant $\Bihol$-cover of $M$ {where each $U_i$ is connected}. Then $\abs{ \mmod{\mathcal{U}}{G}}$ is a cofibrant replacement for the quotient $\infty$-prestack $ \mmod{M}{G}$ in the local projective model structure. Moreover, if $G$ is the trivial group then this reduces to the statement that the \v{C}ech nerve of $\mathcal{U}$ cofibrantly replaces $\varepsilon M$.
\end{proposition}
 
\begin{proof}
Since $\abs{\mmod{\mathcal{U}}{G}}$ is given level-wise as a coproduct of representables and since only certain indices for the coproduct correspond to degenerate summands, then $\abs{\mmod{\mathcal{U}}{G}}$ has \emph{split degeneracies} \cite[definition 4.8]{DHI} and cofibrancy follows. Moreover a morphism in the site never sends non-degenerate summands to degenerate ones for these particular $\infty$-prestacks. Now, by \cite[theorem 6.2]{DHI}, which states the the weak equivalences in the Bousfield localization $\sPre(\Bihol)_{/{\mc S}}$ are equal to the topological weak equivalences of Jardine, it remains to show that the natural map $\abs{\mmod{\mathcal{U}}{G}} \xrightarrow{r} \mmod{M}{G}$ is such a weak equivalence. By lemma \ref{LEM UmmodG 2cosk}, checking the topological weak equivalence at $\pi_0$ and $\pi_1$ suffices, which was already done in lemma \ref{LEM iso hmtpy grps}. 
\end{proof}

\begin{remark}
There is an alternative way to see the topological weak equivalence $\abs{ \mmod{\mathcal{U}}{G}} \to  \mmod{M}{G}$ which we sketch here. The constant simplicial presheaf associated to the simplicial set $BG$ fits into a map between short exact sequences of $\infty$-prestacks,
\begin{equation*}
\begin{tikzcd}
1 \arrow[r] & \Cech\mathcal{U} \arrow[r] \arrow[d] & \abs{\mmod{\mathcal{U}}{G}} \arrow[r] \arrow[d] & BG \arrow[d]\arrow[r] & 1 \\
1 \arrow[r] &  \varepsilon M \arrow[r]& {\mmod{M}{G}} \arrow[r] & BG\arrow[r] & 1,
\end{tikzcd}
\end{equation*}
where the left-most and right-most vertical maps are topological weak equivalences. By comparing the long exact sequences of sheaves of homotopy groups \cite[lemma 1.15]{Ja} induced by the top and bottom rows using the Five Lemma, the middle vertical map must also be a topological weak equivalence.
\end{remark}

\begin{definition}\label{DEF good cover}
An open cover $\mathcal{U}:= \{ U_i \subset M \vert U_i \in \Bihol \}$ of an $n$-dimensional complex manifold $M$ is called a \emph{good open cover} if it is a $\Bihol$-cover of $M$ where each $U_i$ is also connected. If $M$ is moreover a $G$-manifold and the good open cover $\mathcal{U}$ is an equivariant $\Bihol$-cover, then it is called an \emph{equivariant good open cover}.
\end{definition}
Proposition \ref{PROP cofib repl M mmod G} is used to compute the derived mapping space $\mathbb{R} \underline{\sPre}\left(\mmod{M}{G}, \mathcal{F}\right)$ as follows. Given an equivariant good open cover $\mc U$ of $M$, then $\check{\mmod{M}{G}} =\abs{\mmod{\mathcal{U}}{G}} \hookrightarrow \mmod{M}{G}$ provides a cofibrant replacement. When $\mathcal{F} =\mathcal{F}^{\dagger}$ is already fibrant, as in the cases of $\mathcal{F} = \holFormP{k}$  by proposition \ref{PROP holForm stack} and $\mathcal{F} = {\Chart}$  by proposition \ref{PROP chart stack}, the mapping space \eqref{EQ RHom} becomes $\mathbb{R}\text{Hom}(\mmod{M}{G}, \mc F) = \underline{\sPre}(\Bihol) \left(\abs{\mmod{\mathcal{U}}{G}} ,\mc F \right)$. The right hand side can further be computed by the totalization yielding:
\begin{equation}\label{EQ RHom U mmod G}
\mathbb{R}\text{Hom}\left(\mmod{M}{G}, \mathcal{F}\right)\simeq \Tot \left(\mathcal{F}_{\bullet} \left( \abs{\mmod{\mathcal{U}}{G}}^{\bullet} \right)\right) 
\end{equation}
where the co-simplicial-simplicial set $\mathcal{F}_{\bullet} \left( \abs{\mmod{\mathcal{U}}{G}}^{\bullet} \right)$ is defined at each co-simplicial level $k$ by:
\begin{equation}
\mathcal{F}_{\bullet} \left( \abs{\mmod{\mathcal{U}}{G}}^{k} \right):=\prod\limits_{\substack{i_0, \ldots, i_k \\ g_1, \ldots, g_{k-1}}} \mathcal{F}\left( U_{\left(i_0, \ldots, i_k\right) \cdot (g_1 \dots g_{k-1})}\right), \end{equation} 
and each cosimplicial map $d^k$ is induced by the composition of $d^{-1,0}_k$ and $d^{0,-1}_k$ described above. Note that for $k=0$ the $g_i$'s are dropped from the indexing.
\begin{example}\label{EX Tot Chart}
Since $\Chart$ is object-wise $0$-coskeletal and by lemma \ref{LEM: cosk tot} its totalization will also be $0$-coskeletal, then a $p$-simplex in $\Tot \left(\Chart_{\bullet} \left( \abs{\mmod{\mathcal{U}}{G}}^{\bullet} \right)\right) $ is completely determined by its $(p+1)$-many vertices. But each vertex in the totalization is a sequence of $k$-simplices $x^0_0, x^1_1, \ldots $, where each $x^k_k$ is in $\Chart_{k} \left( \abs{\mmod{\mathcal{U}}{G}}^{k} \right)$. Again by the coskeletal nature of $\Chart$, then each $x^k_k$ is completely determined by its vertices which by the coherence condition imposed by face and coface maps, those vertices are uniquely determined by the data in $x^0_0$. And so, a $p$-simplex in $\Tot \left(\Chart_{\bullet} \left( \abs{\mmod{\mathcal{U}}{G}}^{\bullet} \right)\right) $ is $(p+1)$-many atlases $\left\{ U_{i} \xrightarrow{\rho^k_{i}} V^k_{i}\subset \mathbb{C}^n\right\}_{i \in I, k \in \{0,\ldots,p\}}$on $M$ with respect to the $G$-equivariant good cover $\mathcal{U}$ and these atlases have no additional relation to each other.

\end{example}
\begin{example}\label{EX Tot DK forms}
The $p$-simplces in $\Tot \left(\extAnti_{\bullet} \left( \abs{\mmod{\mathcal{U}}{G}}^{\bullet} \right)\right) $ are computed by noting that
\begin{equation}  \Tot \left(\extAnti_{\bullet} \left( \abs{\mmod{\mathcal{U}}{G}}^{\bullet} \right)\right)  =  \Tot \left(\DK \left( \extAnti^{\bullet}\right) \left( \abs{\mmod{\mathcal{U}}{G}}^{\bullet} \right)\right) \stackrel{\text{ lemma \ref{LEM Emilios Tot DK} }}{=}\DK  \left( \left( \tau_{\le 0} \circ \tot\right)  \left(  \extAnti^{\bullet} \left( \abs{\mmod{\mathcal{U}}{G}}^{\bullet}\right) \right)\right).
\end{equation}
Where $\tau_{\le 0}$ is the ``smart truncation'' \eqref{EQ smart trunc}. But in our definition the complex $\extAnti^{\bullet}$ was already truncated (since we took $\delbar$-closed forms in degree $0$).
\end{example}

\section{Applications}
While we attribute the proof of theorem \ref{THM main theorem} to Toledo and Tong, the statement is a generalization to an equality of $\infty$-stacks. In this section, we illustrate the flexibility and improvements that come with this generalization by considering three examples of the domain $\infty$-prestack $\mc F$ in \eqref{EQ BM=Todd spre}: 
\begin{example}[$\mc F = \mmod{V}{G}$]
Let $V \subset \mathbb{C}^n$ be an open Euclidean space and $G \le \restr{\nGpd}{V}$ be a subgroup of the group of biholomorphisms of $V$. Then for the quotient $\infty$-prestack (definition \ref{DEF M mmod G}) $\mc F= \mmod{V}{G}$, we have $\mathbb{R}\text{Hom}( \mmod{V}{G}, \Chart) = \Tot\left( \Chart \left( \mmod{V}{G} \right) \right)$ and a vertex in here (see example \ref{EX Tot Chart}) is uniquely determined by a single biholomorphism $\phi:V \xrightarrow{\sim} U \subset \mathbb{C}^n$. A canonical vertex then is the identity map $id: V \xrightarrow{\sim} V$. At the same time we have that the codomain mapping space in \eqref{EQ BM=Todd spre} is computed as the Dold-Kan of group cohomology: $\mathbb{R}\text{Hom}( \mmod{V}{G}, \holFormP{n}) = \DK \left(  C^{\bullet}\left( G, \holFormP{n}(V)\right)\right)$,
whereby $\Hartmap \circ \BMmap (id: V \to V) = \tau_G$, recovering the group-invariant from \cite{GZ}:
\begin{equation}\label{EQ BM=Todd V mmod G}
\begin{tikzcd}[column sep=huge, row sep=large]
\Tot\left( \Chart \left( \mmod{V}{G} \right) \right)
  \arrow{r}[name=U, label=above:$\left(\Hartmap \circ \BMmap\right)_*$]{} & 
 \DK \left(  C^{\bullet}\left( G, \holFormP{n}(V)\right)\right)\\[5pt]
  \text{on canonical vertex:} \quad  \left(id: V \to V\right) \arrow[mapsto]{r} & \tau_G \in Z^n(G, \holFormP{n}(V),
\end{tikzcd}
\end{equation}
where $C^{\bullet}\left( G, \holFormP{n}(V)\right)$ is truncated in order to apply Dold-Kan as in lemma \ref{EQ Emilios Tot DK}. For any other vertex $\phi: V \xrightarrow{\sim} U \subset \mathbb{C}^n$ in the domain (derived) mapping space  \eqref{EQ BM=Todd spre}, $\Hartmap \circ \BMmap (\phi)$ provides another cocycle in the group cohomology complex for $G$, which is the pullback of $\tau_{G'}$ under the group (complex) isomorphism $G \xrightarrow{\sim} G' = \phi \cdot G \cdot \phi^{-1}$. The image of the unique edge connecting these two domain vertices under $\Hartmap \circ \BMmap$, shows that $\tau_G$ and the conjugate of $\tau_{G'}$ are cohomologous.
\end{example}

\begin{example}[$\mc F =\varepsilon M$]
Let $M$ be a manifold so that $\mc F =\varepsilon M$ is the $\infty$-prestack discussed at the beginning of section \ref{SEC eval on man}, whose cofibrant replacement is given by the nerve of a good open cover $\Cech\mathcal{U}$(see proposition \ref{PROP cofib repl M mmod G}). From here the domain mapping space of  \eqref{EQ BM=Todd spre} can be compute by $\mathbb{R}\text{Hom}(M, \Chart) = \mathbb{R}\text{Hom}(\Cech\mathcal{U}, \Chart) = \Tot \Chart \left( \Cech\mathcal{U} \right)$, where the first equality is via cofibrant replacement since $\Chart$ is already fibrant (i.e., an $\infty$-stack), whereas the second equality utilizes that the homotopy limit induced by computing this mapping space is computed by the associated totalization. Recall as well from example \ref{EX Tot Chart} that $p$-simplices in this particular totalization are given by $(p+1)$ many (independent) atlases on $M$ with respect to $\mc U$. Finally, after applying example \ref{EX Tot DK forms} again to the codomain of \eqref{EQ BM=Todd spre}, we obtain the following interpretation:
\begin{equation}\label{EQ BM=Todd M}
\begin{tikzcd}[column sep=huge]
\mathbb{R}\text{Hom}( \varepsilon M, \Chart)
  \arrow[bend left=15]{r}[name=U, label=above:$\left(\Hartmap \circ \BMmap\right)_*$]{}
  \arrow[bend right=15]{r}[name=D,label=below:${\Cfmap{Todd_n}}_*$]{} &
 \DK \left(  \check{C}^{\bullet}\left( \mc U, \holFormP{n}(V)\right)\right)
  \arrow[shorten <=25pt,shorten >=20pt,equal,to path={(U) -- (D)}]{}\\[10pt]
\text{on vertices:} \quad  \text{atlas} \arrow[mapsto]{r}{} & \text{Todd cocycle},\\[-20pt]
\text{on edges:} \quad  \text{a pair of atlases} \arrow[mapsto]{r}{} & \parbox{4cm}{witness of two atlas' Todd cocycles being cohomologous}
\end{tikzcd}
\end{equation}
where the assignment on vertices is described in \cite{TT76}. Moreover, the equality $\left(\Hartmap \circ \BMmap\right)_* ={\Cfmap{Todd_k}}_*$ offers a complete cochain-level interpretation of Hirzebruch-Riemann-Roch.
\end{example}

\begin{example}[$\mc F = \mmod{M}{G}$]\label{EX application M mmod G}
Let $G$ be a discrete group and $M$ be a $G$-manifold. Then for the quotient $\infty$-prestack (definition \ref{DEF M mmod G}) $\mc F= \mmod{M}{G}$, after cofibrant replacement (proposition \ref{PROP cofib repl M mmod G}) we have $\mathbb{R}\text{Hom}( \mmod{M}{G}, \Chart) = \Tot\left( \Chart \left( \abs{\mmod{\mc U}{G}} \right) \right)$ and a vertex in here is given by a $G$-equivariant atlas. The codomain mapping space in \eqref{EQ BM=Todd spre} is computed as the Dold-Kan of a mix between group and \v{C}ech cohomology: $\mathbb{R}\text{Hom}( \abs{\mmod{\mc U}{G}}, \holFormP{n}) = \DK \left(  \check{C}^{\bullet}\left( G, \holFormP{n}(V)\right)\right)$, and now a generalization of the previous two examples is obtained:
\begin{equation}\label{EQ BM=Todd M mmod G}
\begin{tikzcd}[column sep=huge, row sep=large]
 \Tot\left( \Chart \left( \abs{\mmod{\mc U}{G}} \right) \right)
   \arrow[bend left=15]{r}[name=U, label=above:$\left(\Hartmap \circ \BMmap\right)_*$]{}
  \arrow[bend right=15]{r}[name=D,label=below:${\Cfmap{Todd_n}}_*$]{} &
\DK \left(  \check{C}^{\bullet}\left( G, \holFormP{n}(V)\right)\right)
  \arrow[shorten <=26pt,shorten >=23pt,equal,to path={(U) -- (D)}]{}\\[5pt]
\text{on vertices:} \quad  \text{$G$ equivariant atlas $\phi_{\bullet}$}  \arrow[mapsto]{r} & \tau_{\phi_{\bullet} }\in \check{Z}^n(G, \holFormP{n}(V)),\\[-25pt]
\text{on edges:} \quad  \text{$G$-equivariant atlases $\phi_{\bullet},\phi_{\bullet}'$}  \arrow[mapsto]{r} &{\Cfmap{Todd_n}}_* \left( \phi_{\bullet},\phi_{\bullet}' \right)
\end{tikzcd}
\end{equation}
satisfying $\left(\tau_{\phi_{\bullet}} -  \tau_{\phi_{\bullet}'} = d \left( {\Cfmap{Todd_n}}_* \left( \phi_{\bullet},\phi_{\bullet}' \right) \right) \right)$. Here, the equality $\left(\Hartmap \circ \BMmap\right)_* ={\Cfmap{Todd_n}}_*$ offers a complete cochain-level interpretation of an equivariant Hirzebruch-Riemann-Roch. Below, the new results relevant to this example are recorded.
\end{example}

\begin{definition}
A manifold $M$ is an \emph{affine manifold} if there exists an atlas where each transition map is affine; such an atlas is called an \emph{affine atlas}. 
\end{definition}
\begin{corollary}
Let $M$ be a $G$-manifold, $\mc U$ be a $G$-equivariant cover, and $\phi_{\bullet}$ be an affine atlas which $G$ acts on by affine biholomorphisms, then $\tau_{\phi_{\bullet}} = 0$ (see example \ref{EX application M mmod G}).
\begin{proof}
Note that our $\tau_{\phi_{\bullet}}$ is precisely $\omega^n$ from \cite{TT76}. By \cite[theorem 10.9 and lemma 11.3]{TT76} $\omega^n = \tau^{n,0}$ is trivial when evaluated on $(\phi_1, \ldots, \phi_n)$ where one of the $\phi_i$ is affine. Since the conditions of the corollary will always force all of the $\phi_i$ to be affine then the statement follows.
\end{proof}
\end{corollary}
Given the above corollary, the next one is an immediate consequence of the fact (see example \ref{EX application M mmod G}) that two $G$-equivariant atlases with respect to the same cover result in cohomologous invariants.
\begin{corollary}
The invariant $\tau$ from example \ref{EX application M mmod G} is an obstruction to finding an affine atlas which the groups acts affinely on.
\end{corollary}
\appendix

\end{document}